\documentclass[smallextended]{svjour3} 
\smartqed 
\usepackage{amssymb}
\usepackage{amsmath}
\usepackage{epsfig}
\usepackage{graphbox}
\usepackage{color}
\usepackage{mathrsfs}
\usepackage[sort&compress,numbers]{natbib}
\usepackage{mathrsfs}
\usepackage{ulem}
\usepackage{enumerate}
\usepackage{enumitem}
\usepackage[numbers]{natbib}
\usepackage{lmodern}

\usepackage{txfonts}
\usepackage{latexsym}

\newcommand*{\vv}[1]{\vec{\mkern0mu#1}}
\newcommand{\norm}[1]{\Vert#1\Vert}
\newcommand{\bR}{{\mathbb R}}
\newcommand{\bN}{{\mathbb N}}
\newcommand{\bI}{\mathbb{I}}
\newcommand{\bX}{V_{\partial}}
\newcommand{\bU}{\mathbb{U}}
\newcommand{\bP}{\mathbb{P}}
\newcommand{\tD}{\mathbb{D}}

\newcommand{\StabV}{{\tt StabV}}
\newcommand{\StabN}{{\tt Stab}}
\newcommand{\EquiL}{{\tt Equid}}
\newcommand{\EquiV}{{\tt EquidV}}

\newcommand{\mX}{\mathscr{X}}
\newcommand{\mZ}{{\mathfrak{x}}}
\newcommand{\mR}{\mathcal{R}}
\newcommand{\mS}{S}

\newcommand{\bG}{\nabla_c}
\newcommand{\bD}{\mat{\tD_c}}
\newcommand{\dV}{\mathrm{d}V}
\newcommand{\drz}{\mathrm{d}r\mathrm{d}z}
\newcommand{\vol}{\operatorname{vol}}

\newcommand{\pp}[2]{\frac{\partial#2}{\partial#1}}

\newcommand{\rd}{\;{\rm d}}
\newcommand{\id}{{\rm id}}
\newcommand{\dd}[1]{\frac{\rm d}{{\rm d}#1}}
\newcommand{\ddt}{\dd{t}}
\newcommand{\nn}{\nonumber}
\newcommand{\ttau}{\Delta t}

\newcommand{\cira}{\mbox{$s\!\!\!\!\:/$}}

\newcommand{\mat}[1]{\uuline{#1}}
\newcommand{\adapt}[2]{{\rm adapt}_{#1,#2}}

\journalname{J. Sci. Comput.}

\begin{document}
\title{Unfitted finite element methods for axisymmetric two-phase flow}
\titlerunning{Unfitted FEMs for axisymmetric two-phase flow}  

\author{Harald Garcke \and 
       Robert N\"urnberg  \and 
       Quan~Zhao
}
\authorrunning{H. Garcke, R. N\"urnberg and Q. Zhao}

\institute{Harald Garcke \at
              Fakult{\"a}t f{\"u}r Mathematik,
              Universit{\"a}t Regensburg, 
              93040 Regensburg, Germany.\\
              \email{harald.garcke@ur.de}          
           \and
           Robert N\"urnberg \at
            Dipartimento di Mathematica, 
            Universit\`a di Trento,
            38123 Trento, Italy.\\
            \email{robert.nurnberg@unitn.it}
            \and 
            Quan Zhao \at 
            Fakult{\"a}t f{\"u}r Mathematik,
            Universit{\"a}t Regensburg, 
            93040 Regensburg, Germany.\\
              \email{quan.zhao@ur.de} 
}

\date{Received: date / Accepted: date}

\maketitle

\begin{abstract}
We propose and analyze unfitted finite element approximations for the two-phase incompressible Navier--Stokes flow in an axisymmetric setting. The discretized schemes are based on an Eulerian weak formulation for the Navier--Stokes equation in the 2d-meridian halfplane, together with a parametric formulation for the generating curve of the evolving interface. We use lowest order Taylor--Hood and piecewise linear elements for discretizing the Navier--Stokes formulation in the bulk and the moving interface, respectively. 
We discuss a variety of schemes, amongst which is a linear scheme that enjoys
an equidistribution property on the discrete interface and good volume
conservation. An alternative scheme can be shown to be unconditionally stable
and to conserve the volume of the two phases exactly. Numerical results are presented to show the robustness and accuracy of the introduced methods for simulating both rising bubble and oscillating droplet experiments.
\end{abstract}

\renewcommand{\thefootnote}{\arabic{footnote}}
\setlength{\parindent}{2em}

\setcounter{equation}{0}
\section{Introduction} \label{sec:intro}
\setlength\parindent{24pt}

Two-phase flows are widely observed phenomena in nature, and they have important applications in industrial engineering and scientific experiments. Numerical approximations of these flows have attracted a lot of attention, and a large body of numerical methods have been introduced in recent decades, see e.g., \cite{Hirt1981volume,Hughes81,UnverdiT92,Sussman94level,Anderson1998,Sethian99level,Tryg01,Bansch01,osher02level,Renardy02,Perot2003moving,Feng2006fully,Ganesan06,Gross07extended,olsson07,Styles2008finite,Popinet09,BGN2013eliminating,Grun2014two,BGN15stable,GarckeHK16,Cheung18mass, FrachonZ19,Agnese20,Zhao2020energy,Duan2022energy, GNZ23}. 
Nevertheless, accurate and efficient numerical approximations for fully three-dimensional flows remain a significant challenge. This is despite several 
different treatments of the moving interface being available in the literature. Among these one of the prominent methods is the so-called level set method \cite{Sussman94level,Sethian99level,osher02level,olsson07,Gross07extended}, where the interface is represented as a level set function. The main drawback in this approach is that the level set function needs to be re-initialized constantly in order to stay close to the signed distance function. This usually introduces unphysical shifting of the interface, thus causing significant volume loss \cite{Gibou2018review}. 
On the other hand, the front-tracking method, which tracks the interface explicitly with a set of markers or a lower-dimensional mesh, provides a very accurate approximation of the interface and in particular the volume preservation aspect, see \cite{BGN2013eliminating,BGN15stable,GNZ23}. However, the front tracking method may struggle with complicated mesh manipulations which are often necessary, for example, in the computation of interface-bulk cross terms in the unfitted mesh approach \cite{BGN10stefan,BGN15stable} or during the dynamic control of the bulk and interface meshes in the moving fitted mesh approach \cite{Quan07moving,Anjos20143d}.  
 
The computationally challenging aspect of fully three-dimensional two-phase
flows not only comes from the large size of the problem itself but also involves the possible difficulties of the numerical methods themselves, see e.g., \cite{Mirjalili2017, GrossR11, Gibou2018review}. Fortunately, in many situations the two-phase flow satisfies a rotational symmetry, meaning that the complex 3d problem can be reduced to a much simpler two-dimensional problem in the meridian halfplane. Moreover, the evolving fluid interface can be modelled by a one-dimensional generating curve, dramatically reducing the computational complexity. For front tracking methods, in particular, dealing with a discrete curve avoids undesirable mesh distortions as well as the associated remedies like mesh quality controls. There exist several numerical approximations for the axisymmetric two-phase incompressible flows, including the level set method \cite{Sussman00coupled,Chessa2003}, the diffuse-interface approach \cite{Kim05diffuse,Huang22diffuse} and the moving mesh interface-tracking method \cite{Ganesan08accurate,Gros18interface}. The main aim of this work is to explore unfitted finite element approximations for the axisymmetric two-phase incompressible flow and pay particular attention to energy stability and volume preservation aspects. The numerical methods we consider fall into the framework of a parametric approach, which was introduced in \cite{BGN2013eliminating,BGN15stable} by Barrett, Garcke and N\"urnberg, see \eqref{eq:3dwf3}-\eqref{eq:3dwf4}.
In this `BGN' approach, the surface tension force is accurately computed with the help of a curvature identity, see \cite{Dziuk90}. For the evolving interface, the normal velocity is prescribed by the normal part of the fluid velocity, which allows tangential degrees of freedoms to improve the quality of the interface mesh. In particular, for an evolving curve, a semidiscrete version of the discretized scheme leads to uniformly distributed mesh points \cite{BGN07}. 
 The interested reader is referred to the recent review paper \cite{Barrett20} for a detailed discussion of this tangential velocity. We also note that, in the context of geometric flows, alternative approaches that also introduce a benign tangential velocity have been considered in \cite{m17,Hu22evolving}.

The unfitted finite element approximations we introduce in this paper crucially
rely on an appropriate formulation of the surface tension forces in the
axisymmetric setting. Here we rely on two possible formulations for the mean
curvature of the interface presented in \cite{BGN19asy,BGN19variational}.
The first formulation leads to a discretization of the curvature of the
generating curve, thus providing an equidistribution property as in the
two-dimensional case, see \cite{BGN07,BGN15stable}. The alternative is to
approximate the mean curvature of the three-dimensional interface, which has
the advantage that an unconditional stability estimate can be shown.
The overall numerical methods are then obtained by coupling the interface
evolution to the fluid flow in the bulk. We follow the approach from 
the unfitted 2d/3d approximations of the Eulerian weak formation that was 
introduced in \cite{BGN15stable}. In particular, by employing a special
XFEM enrichment procedure for the pressure space, a volume preservation
property can be shown on the semidiscrete level. 

Strengthening the volume conservation property to the fully discrete level is
also possible, on combining the XFEM technique with discrete time-weighted 
interface normals as derived in \cite{BGNZ22volume}.
These time-weighted interface approximations were initially explored in the context of the surface diffusion flow \cite{Jiang21,BZ21SPFEM} and then extended to the axisymmetric geometric flows in \cite{BGNZ22volume}. In the 2d/3d setting 
an exactly volume-preserving finite element method for the two-phase incompressible flow was recently proposed by the authors in \cite{GNZ23}. 

The rest of the paper is organized as follows. In Section~\ref{sec:3df}, we review an Eulerian weak formulation for the sharp-interface model of two-phase incompressible flow in three dimensions. Next, in Section~\ref{sec:25df}, we introduce two axisymmetric weak formulations for the considered flow under a suitable axisymmetric setting. Subsequently, in Section~\ref{sec:uffem}, we propose four discretized schemes based on the axisymmetric formulations and analyze their properties of equidistribution, energy stability and volume preservation. Numerical results for rising bubble and oscillating droplet experiments are presented in Section~\ref{sec:num} to show the robustness and the accuracy of the introduced schemes. Finally, we draw some conclusions in Section~\ref{sec:con}.

\setcounter{equation}{0}
\section{Mathematical formulations in 3d}\label{sec:3df}

We consider the dynamics of the two-phase fluids in a fixed bounded domain $\Omega\subset\bR^3$ with  
$\Omega=\Omega_+(t)\cup\Omega_-(t)\cup\mS(t)$, where $\Omega_\pm(t)\subset\bR^3$ are the regions occupied by the two fluids, $t$ is the time variable, and $\mS(t)$ is a closed surface representing the fluid interface which separates the two regions. 
In the following we stay close to the presentation in \cite{BGN15stable},
see also \cite{GrossR11}.  
We assume a parameterization of $\mS(t)$ is given by: $\vec\Upsilon(\cdot,t):\mathcal{O}\to\bR^3$, where $\mathcal{O}$ is the reference surface. The induced interface velocity is then defined as
\begin{equation}
\mathcal{\vv V}(\vec\Upsilon(\vec q,t), t) = \partial_t\vec\Upsilon(\vec q,t)\quad\forall\vec q\in\mathcal{O}.\nn
\end{equation}
Moreover, we introduce the mean curvature of the interface as
\begin{equation}
\mathcal{H}\,\vec n_{_\mS} = \Delta_s\vec\id,
\end{equation}
where $\vec n_{_\mS}$ is the unit normal to $\mS(t)$ and points to the region $\Omega_+(t)$, $\Delta_s=\nabla_s\cdot\nabla_s$ is the surface Laplacian operator on the surface with $\nabla_s$ being the surface gradient, and $\vec\id$ is the identity function on $\mS(t)$.

Let $\breve{\vec u}(\cdot, t):\Omega\times[0,~T]\to\mathbb{R}^3$ be the fluid velocity,  and $\breve{p}(\cdot, t):\Omega\times[0, T]\to\mathbb{R}$ be the pressure. Moreover, we use $\rho_\pm$ and $\mu_\pm$ to denote the densities and viscosities of the fluids in $\Omega_\pm(t)$, respectively. Then the governing equations for the two-phase Navier--Stokes flow  are given by 
\begin{subequations}
\label{eqn:fluidy}
\begin{alignat}{3}
\label{eq:fluid1}
\rho_\pm\bigl(\partial_t\breve{\vec u} + \breve{\vec u}\cdot\nabla\breve{\vec u}\bigr) 
& = \nabla\cdot\mat{\sigma} + \rho_\pm\,\breve{\vec g}\quad &&\mbox{in}\;\Omega_\pm(t),\\ 
\label{eq:fluid2}
 \nabla\cdot\breve{\vec u} & = 0\quad &&\mbox{in}\;\Omega_\pm(t),\\
 \label{eq:fluid3}
 [\breve{\vec u}]_-^+ &= 0\quad &&\mbox{on}\;\mS(t),\\
  \label{eq:fluid4}
  \bigl[\mat{\sigma}\,\vec n_{_\mS}\bigr]_-^+&= -\gamma\mathcal{H}\,\vec n_{_\mS} &&\mbox{on}\;\mS(t),\\
  \label{eq:fluid5}
 \mathcal{\vv V}\cdot\vec n_{_\mS}&= \breve{\vec u}\cdot\vec n_{_\mS} &&\mbox{on}\;\mS(t),\\
  \label{eq:fluid6}
 \breve{\vec u} &= \vec 0\quad &&\mbox{on}\;\partial_{_1}\Omega,\\
\breve{\vec u}\cdot\vec n_{_b}=0 ,\quad(\mat{\sigma}\,\vec n_{_b})\cdot\vec t&=0\quad\forall\vec t\in\{\vec n_{_b}\}^\perp\quad&&\mbox{on}\;\partial_{_2}\Omega.
 \label{eq:fluid7}
\end{alignat}
\end{subequations}
Here \eqref{eq:fluid1}-\eqref{eq:fluid2} are the two-phase Navier--Stokes equations, $\breve{\vec g}$ is the body acceleration, $\mat{\sigma}$ is the stress tensor defined by 
\begin{equation}
\mat{\sigma}= \mu_\pm[\nabla\breve{\vec u} + (\nabla\breve{\vec u})^T]-\breve{p}\,\mat{I_3}=2\mu_\pm\,\mat\tD(\breve{\vec u}) - \breve{p}\,\mat{I_3}\quad\mbox{in}\quad\Omega_\pm(t),\nn
\end{equation}
with $\mat\tD(\breve{\vec u})=\frac{1}{2}[\nabla\breve{\vec u} + (\nabla\breve{\vec u})^T]$ being the strain rate,  and $\mat{I_3}\in\mathbb{R}^{3\times 3}$ is the identity matrix. Equations \eqref{eq:fluid3}--\eqref{eq:fluid5} are the classical interface conditions, where $[\cdot]_-^+$ denotes the jump value from $\Omega_-(t)$ to $\Omega_+(t)$, $\gamma$ is the surface tension of $\mS(t)$. Moreover, we prescribe the no slip condition \eqref{eq:fluid6} on $\partial_1\Omega$ and a free slip condition \eqref{eq:fluid7} on $\partial_2\Omega$, where $\vec n_b$ is the outer unit normal to $\partial\Omega$, and $\{\vec n_b\}^\perp:=\bigl\{\vec t\in\bR^3\;:\;\vec t\cdot\vec n_b=0\bigr\}$. 

We define the total energy of the system as
\begin{equation}
\mathcal{E}(t)
 = \tfrac{1}{2}\int_{\Omega_+(t)}\rho_+\,|\breve{\vec u}|^2\,\dV+\tfrac{1}{2}\int_{\Omega_-(t)}\rho_-\,|\breve{\vec u}|^2\,\dV+\gamma|\mS(t)|,\label{eq:3DEnergy}
\end{equation}
where $|\mS(t)| = \int_{\mS(t)}1\,\rd S$.
Then the dynamic system obeys the volume conservation and energy laws as follows
\begin{subequations}\label{eqn:3dlaws}
\begin{align}
&\ddt\vol(\Omega_-(t))=\int_{\mS(t)}\mathcal{\vv V}\cdot\vec n_{_\mS}\,\rd S=\int_{\Omega_-(t)}\nabla\cdot\breve{\vec u}\,\dV = 0,\\
&\ddt\mathcal{E}(t) = -2\int_\Omega\mu\,\mat{\tD}(\breve{\vec u}):\mat{\tD}(\breve{\vec u})\,\dV + \int_\Omega\rho\,\breve{\vec g}\cdot\breve{\vec u}\,\dV,
\end{align}
\end{subequations}
where we introduced $\rho(\cdot, t) = \rho_+\mX_{_{\Omega_+(t)}} + \rho_-\mX_{_{\Omega_-(t)}}$ and $\mu(\cdot, t) = \mu_+\mX_{_{\Omega_+(t)}} + \mu_-\mX_{_{\Omega_-(t)}}$ 
with  $\mX_{_E}$ being the characteristic function of a set $E$.

We denote by $(\cdot,\cdot)_{\Omega}$ and $\langle\cdot,\cdot\rangle_{\mS(t)}$ the $L^2$-inner products over $\Omega$ and $\mS(t)$, respectively. On defining the function spaces 
\begin{subequations}
\begin{align}
\breve{\bU} &= \left\{\vec\chi\in [H^1(\Omega)]^3\,:\,\vec\chi = 0\;\mbox{on}\;\partial_1\Omega,\quad \vec\chi\cdot\vec n_b = 0\;\mbox{on}\;\partial_2\Omega\right\},\nn\\
\breve{\bP} &= \bigl\{\eta\in L^2(\Omega)\,:\,(\eta,~1)_{\Omega} =0\bigr\},\quad\breve{\mathbb{V}}:=H^1(0,T;[L^2(\Omega)]^3)\cap L^2(0,T;\breve{\bU}),\nn
\end{align}
\end{subequations}
we then recall the following Eulerian weak formulation for the two-phase incompressible flow from \cite[(3.9)]{BGN15stable}. Given the initial velocity $\breve{\vec u}_0 = \breve{\vec u}(\cdot,0)$ and the fluid interface $\mS(0)$, we find $\mS(t)=\vec\Upsilon(\mathcal{O},t)$ with $\mathcal{\vv V}(\cdot, t)\in [L^2(\mS(t))]^3$, $\breve{\vec u}\in \breve{\mathbb{V}}$, $\breve{p}\in L^2(0,T;\breve{\bP})$, and $\mathcal{H}(\cdot, t)\in L^2(\mS(t))$ such that for all $t\in(0,T]$
 \begin{subequations}\label{eqn:3dwf}
\begin{alignat}{3}
\label{eq:3dwf1}
&\tfrac{1}{2}\left[\ddt\bigl(\rho\,\breve{\vec u},~\breve{\vec\chi}\bigr)_\Omega + \bigl(\rho\,\partial_t\breve{\vec u},~\breve{\vec\chi}\bigr)_\Omega - \bigl(\rho\breve{\vec u},~\partial_t\breve{\vec\chi}\bigr)_\Omega\right] + 2\bigl(\mu\,\mat\tD(\breve{\vec u}),~\mat\tD(\breve{\vec\chi})\bigr)_\Omega\nn\\
&\hspace{1cm}+ \tfrac{1}{2}\left[\bigl(\rho\,[\breve{\vec u}\cdot\nabla]\breve{\vec u},~\breve{\vec\chi}\bigr)_\Omega-\bigl(\rho\,[\breve{\vec u}\cdot\nabla]\breve{\vec\chi},~\breve{\vec u}\bigr)_\Omega\right] -\bigl(\breve{p},~\nabla\cdot\breve{\vec\chi}\bigr)_\Omega \nn\\
&\hspace{1cm}- \gamma\langle\mathcal{H}\,\vec n_{_\mS},~\breve{\vec\chi}\rangle_{\mS(t)} = \bigl(\rho\,\breve{\vec g},~\breve{\vec\chi}\bigr)_\Omega\qquad\forall\breve{\vec\chi}\in\breve{\mathbb{V}},\\[0.2em]
\label{eq:3dwf2}
&\hspace{0.5cm}\bigl(\nabla\cdot\breve{\vec u},~\breve{q}\bigr)_\Omega=0\qquad\forall \breve{q}\in\breve{\bP},\\[0.2em]
\label{eq:3dwf3}
&\hspace{0.5cm}\big\langle\mathcal{\vv V}-\breve{\vec u},~\breve{\zeta}\,\vec n_{_\mS}\big\rangle_{\mS(t)}=0\qquad\forall\breve{\zeta}\in L^2(\mS(t)),\\[0.2em]
&\hspace{0.5cm}\big\langle\mathcal{H}\,\vec n_{_\mS},~\breve{\vec\eta}\big\rangle_{\mS(t)} + \big\langle\nabla_s\vec\id,~\nabla_s\breve{\vec\eta}\big\rangle_{\mS(t)}=0\qquad\forall\breve{\vec\eta}\in [H^1(\mS(t))]^3.\label{eq:3dwf4}
\end{alignat}
\end{subequations}

\setcounter{equation}{0}
\section{Axisymmetric formulations}\label{sec:25df}

\subsection{The axisymmetric setting}\label{sec:asyset}
\begin{figure}[!htp]
\centering
\includegraphics[width=0.8\textwidth]{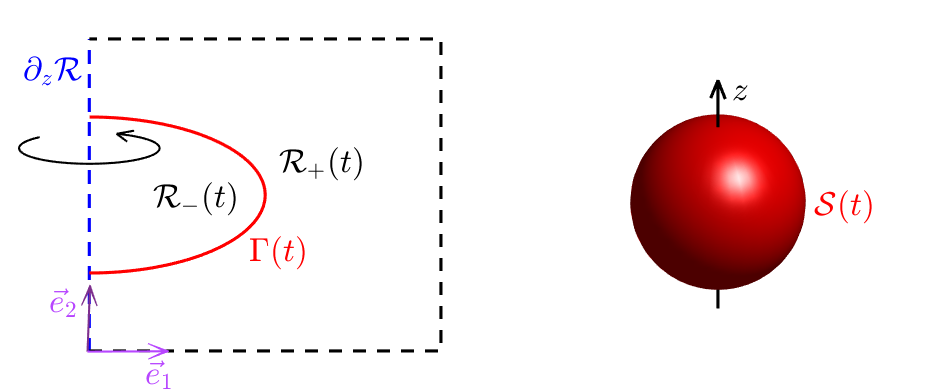}
\caption{Left panel: A schematic illustration of the meridian halfplane $\mR=\mR_-(t)\cup\mR_+(t)\cup\Gamma(t)$. Right panel: The axisymmetric surface $\mS(t)$.}
\label{fig:asyd}
\end{figure}

We assume that the considered flow is rotationally symmetric with respect to the $z$-axis, as shown in Fig.~\ref{fig:asyd}. We then start from the 3d formulation \eqref{eqn:3dwf} and derive the axisymmetric weak formulation in  the cylindrical coordinates. In this setting,  we have the following three assumptions.
\begin{enumerate}[label=$\mathbf{A \arabic*}$, ref = $\mathbf{A \arabic*}$]
\item \label{asp:I} The bounded domain $\Omega$ can be obtained by rotating the 
subset $\mR$ of the meridian halfplane with respect to the $z$-axis. Here
$\mR=\mR_+(t)\cup\mR_-(t)\cup\Gamma(t)$, where $\mR_\pm(t)$ are the corresponding rotated sets for $\Omega_\pm(t)$, and $\Gamma(t)$ is the generating curve of the axisymmetric interface $\mS(t)$.
\item \label{asp:II} We have a natural decomposition of the boundary of $\mR$ as $\partial\mR = \partial_1\mR\cup\partial_2\mR\cup\partial_z\mR$, 
where $\partial_i\mR$ corresponds to $\partial_i\Omega$ for $i=1,2$, and $\partial_z\mR$ is the artificial boundary of $\mR$ on the $z$-axis. Besides, $\Gamma(t)$ is assumed to be  a smooth open curve whose two end points are attached to $\partial_z\mR$ with a $90^\circ$ contact angle condition, meaning that $\mS(t)$ is a  smooth genus-0 surface without boundary.
\item \label{asp:III} There is no rotational fluid velocity and the velocity components in the $r,z$-directions are independent of the azimuthal angle $\varphi$. Besides, the pressure is independent of $\varphi$. 
\end{enumerate}

By \ref{asp:I} and using the cylindrical coordinates, it is natural to consider the  following parameterization of the interface $\mS(t)$
\begin{equation}\label{eq:mSp}
\vec\Upsilon: (\alpha,\varphi, t)\mapsto\bigl(\mZ^r(\alpha, t)\cos\varphi,  \mZ^z(\alpha, t),\mZ^r(\alpha, t)\sin\varphi\bigr)^T,\quad\alpha\in\bI,\quad \varphi\in[0,2\pi],
\end{equation}
for $ t\in[0,T]$,  where $\bI = (0,1)$ with $\partial\bI =\{0,1\}$. Then we could extract the $r,z$-components to form a parameterization of the generating curve $\Gamma(t)$ as
\begin{equation}
\vec\mZ(\alpha, t)= \bigl(\mZ^r(\alpha, t), \mZ^z(\alpha, t)\bigr)^T: \overline{\bI}\times[0,T]\to\bR_{\geq 0}\times \bR.
\label{eq:gammap}
\end{equation}
By \ref{asp:II}, the attachment of the two end points on the $z$-axis implies that 
\begin{equation}
\vec\mZ(\alpha,t)\cdot\vec e_1=0\quad\forall\alpha\in\partial\bI.
\label{eq:gammabd}
\end{equation}
We further assume that $|\vec\mZ_\alpha(\cdot,t)|>0$ and $\vec\mZ(\cdot,t)\cdot\vec e_1>0$ in $\bI$. On introducing the arc length $s$ of the curve $\Gamma(t)$, i.e., $\partial_s =
|\vec{\mZ}_\alpha|^{-1}\,\partial_\alpha$, the unit tangent and unit
normal to the curve $\Gamma(t)$ are defined as
\begin{equation}
\vec\tau = \vec\mZ_s,\qquad \vec\nu = -\vec\tau^\perp,\nn
\end{equation}
where $\cdot^\perp$ denotes a clockwise rotation by $\frac{\pi}{2}$, and $\vec\nu$ is the unit normal to $\Gamma$ and points into $\mR_+(t)$. We introduce the mean curvature of  $\mS(t)$ , 
\begin{equation}
\varkappa = \kappa - \frac{\vec\nu\cdot\vec e_1}{\vec\mZ\cdot\vec e_1}\quad\mbox{with}\quad \kappa\,\vec\nu = \vec\mZ_{ss}, \label{eq:kappaEQ}
\end{equation}
where $\kappa$ is the curvature of the generating curve $\Gamma(t)$. Clearly, for $\varkappa$ to remain bounded we need that 
$\vec\nu \cdot \vec e_1 = 0$ on $\partial\bI$, which is equivalent to
\begin{equation} \label{eq:bc}
\vec\tau \cdot \vec e_2 = 0
\quad \text{on } \partial\bI,
\end{equation}
i.e.\ the $90^\circ$ contact angle condition from \ref{asp:II}.

Furthermore, by \ref{asp:III} we can write the fluid velocity in terms of the cylindrical coordinates as
\begin{equation} \label{eq:uu}
\breve{\vec u}(\vec x, t)= \bigl(u^r(r,z,t)\cos\varphi, u^z(r,z,t), u^r(r,z,t)\sin\varphi\bigr)^T.
\end{equation}
Then we could extract the $r,z$ components to form a vector field in $\mR$ as
\begin{equation}
\vec u(r,z,t) = \bigl(u^r(r,z,t), u^z(r,z,t)\bigr)^T:\mR\times[0,T]\to\bR^2.
\end{equation}
We further introduce
\begin{equation}
\nabla_c = \vec e_1\,\partial_r+ \vec e_2\,\partial_z, \qquad \mat{\tD_c}(\vec u) = \tfrac{1}{2}[\nabla_c\vec u + (\nabla_c\vec u)^T].\label{eq:nbc}
\end{equation}

\subsection{Axisymmetric weak formulations}

By \eqref{eq:gammabd}, we introduce the function space for $\vec\mZ(\cdot,t)$
as
\begin{equation}
\bX=\Bigl\{\vec\eta \in [H^1(\bI)]^2 :\; \vec\eta(\alpha)\cdot\vec e_1 = 0
\;\; \forall\ \alpha \in \partial\bI\Bigr\}.
\end{equation}
Moreover, we define the weighted Sobolev spaces
\begin{align}
&L_a^2(\mR) := \Bigl\{\chi: \int_{\mR}r^a\,\chi^2\drz  < +\infty\Bigr\},\quad a\in\bR,\nn\\
&H_a^1(\mR):=\Bigl\{\chi\in L_a^2(\mR): \partial_r\chi\in L_a^2(\mR), \;\partial_z \chi\in L_a^2(\mR)\Bigr\},\nn 
\end{align}
and then introduce the function spaces for the velocity and pressure (see~\cite{Bernardi99,Belhachmi06weighted})
\begin{subequations}
\begin{align}
\bU &:= \left\{\vec\chi\in [H_1^1(\mR)]^2:\;(\vec\chi\cdot\vec e_1)\in L_{-1}^2(\mR),\;\vec\chi=\vec 0\;\;\mbox{on}\;\;\partial_{1}\mR,\;\;\vec \chi\cdot\vec n = 0\;\;\mbox{on}\;\;\partial_{2}\mR\right\},\nn\\
\mathbb{V}&:=H^1(0,T; [L_1^2(\mR)]^2)\cap L^2(0,T;\bU),\quad 
\bP := \Bigl\{\chi\in L^2_1(\mR):  \int_\mR r\,\chi\,\drz=0\Bigr\},\nn
\end{align}
\end{subequations}
where $\vec n$ is the outer unit normal to $\partial\mR$.
Let $(\cdot,\cdot)$ and $\langle\cdot,\cdot\rangle$ be the $L^2$-inner products over $\mR$ and $\bI$, respectively. Then by imposing the axisymmetric conditions, we are able to transform the 3d formulation \eqref{eqn:3dwf} into a set of equations in the subset $\mR$ of the 2d meridian halfplane and in $\bI$ as follows. Given the initial velocity $\vec u_0\in\bU$ and $\vec\mZ(\cdot,0)\in\bX$, we find $\vec u\in\mathbb{V}$, $p\in L^2(0,T;\bP)$, $\vec\mZ\in H^1(0,T;\bX)$ and $\varkappa\in L^2(0,T; L^2(\bI))$ such that for all $t\in(0,T]$
\begin{subequations}\label{eqn:weakform}
\begin{align}
&\tfrac{1}{2}\left[\ddt\bigl(\rho_{_c}\vec u,~\vec\chi\,r\bigr) + \bigl(\rho_{_c}\partial_t\vec u,~\vec\chi\,r\bigr) - \bigl(\rho_{_c}\vec u,~\partial_t\vec\chi\,r\bigr)\right]+\mathscr{A}(\rho_{_c},\vec u;\vec u, \vec\chi) \nn\\
&\hspace{1cm} + 2\bigl(\mu_{_c}\,r^{-1}\,[\vec u\cdot\vec e_1],~\vec\chi\cdot\vec e_1\bigr)+ 2\bigl(\mu_{_c}\,r\,\mat{\tD_c}(\vec u),~\mat{\tD_c}(\vec\chi)\bigr)- \bigl(p,~\bG\cdot[r\,\vec\chi]\bigr)\nn\\
&\hspace{1cm} -\,\gamma\big\langle\,[\vec\mZ\cdot\vec e_1]\,\varkappa\,\vec\nu,~\vec\chi\,|\vec\mZ_\alpha|\,\big\rangle = \bigl(\rho_{_c}r\,\vec g,~\vec\chi\bigr)\qquad\forall\vec\chi\in\mathbb{V},\label{eq:weak1}\\[0.5em]
\label{eq:weak2}
&\hspace{0.1cm}\bigl(\bG\cdot[r\,\vec u], ~q\bigr)=0\qquad\forall q\in\bP,\\[0.5em]
\label{eq:weak3}
&\hspace{0.1cm}\big\langle(\vec\mZ\cdot\vec e_1)\,[\partial_t \vec\mZ- \vec u],~\zeta\,\vec\nu\,|\vec\mZ_\alpha|\big\rangle = 0\qquad\forall \zeta\in L^2(\bI),\\[0.5em]
\label{eq:weak4}
&\hspace{0.1cm}\big\langle(\vec\mZ\cdot\vec e_1)\,\varkappa\,\vec\nu,~\vec\eta\,|\vec\mZ_\alpha|\big\rangle+\big\langle\vec e_1\cdot\vec\eta,~|\vec\mZ_\alpha|\big\rangle + \Big\langle\frac{(\vec\mZ\cdot\vec e_1)}{|\vec\mZ_\alpha|}\,\vec\mZ_\alpha,~\vec\eta_\alpha\Big\rangle=0\qquad\forall\vec\eta\in\bX,
\end{align}
\end{subequations}
where in these expressions we introduced
\begin{subequations}
 \begin{align}
 \label{eq:rhoc}
 &\rho_{_c}(\cdot,t)  = \rho_+\mX_{_{\mR_+(t)}}+ \rho_-\mX_{_{\mR_-(t)}},\quad \mu_{_c}(\cdot,t)  = \mu_+\mX_{_{\mR_+(t)}}+ \mu_-\mX_{_{\mR_-(t)}},\\
 &\mathscr{A}(\rho_{_c},\vec v; \vec u, \vec \chi) = \tfrac{1}{2}\bigl(\rho_{_c}\,r,~[\vec v\cdot\bG]\vec u\cdot\vec\chi-[\vec v\cdot\bG]\vec\chi\cdot\vec u\bigr).
 \label{eq:asyt}
\end{align}
\end{subequations}
Compared to \eqref{eqn:3dwf}, we note an extra weight $r$ or $\vec\mZ\cdot\vec e_1$ was introduced in \eqref{eqn:weakform} since the volume and surface integrals were transformed into area and line integrals. The detailed derivation of \eqref{eqn:weakform} is presented in Appendix~\ref{app:asywf}. The main advantage of this axisymmetric reduction is that the boundary conditions at the artificial boundary appear naturally in the formulation,  see \cite{Ganesan08accurate}. For example, we actually have the free slip boundary condition for the velocity on $\partial_z\mR$, which is inferred by $\vec u\in\mathbb{V}$. Moreover, it can be shown that \eqref{eq:weak4} implies the $90^\circ$ contact angle condition \eqref{eq:bc}, see \cite[Appendix~A]{BGN19variational} for details.

In this axisymmetric setting, we denote by $E(\rho_{_c}, \vec u(t),\vec\mZ(t))$ and $M(\vec\mZ(t))$ the energy of the system and the enclosed volume of the axisymmetric surface $\mS(t)$, respectively.  On recalling \eqref{eq:3DEnergy} and using \eqref{eq:mSp} and \eqref{eq:uu}, then it holds
\begin{subequations}
\begin{align}
\label{eq:asyE}
E(\rho_{_c}, \vec u(t), \vec\mZ(t)) &= \pi\Bigl(\rho_{_c}|\vec u|^2,~ r\Bigr) 
+ 2\pi\gamma\Big\langle\vec\mZ\cdot\vec e_1,|\vec\mZ_\alpha|\Big\rangle,\\
\label{eq:asyM}
M(\vec\mZ(t))&=\pi\left\langle (\vec\mZ \cdot\vec e_1)^2\, 
\vec\nu, \vec e_1\,|\vec\mZ_\alpha|\right\rangle,
\end{align}
\end{subequations}
where the reader can refer to \cite[(3.10)]{BGN19asy} for \eqref{eq:asyM}. On recalling \eqref{eqn:3dlaws} and using direct calculation it is not difficult to show that solutions to \eqref{eqn:weakform} satisfy
\begin{subequations}
\begin{align}
\label{eq:dtE}
\ddt E(\vec u(t),\vec\mZ(t)) &= \pi\ddt\Bigl(\rho_{_c}|\vec u|^2,~ r\Bigr) + 2\pi\gamma\left\langle\vec\mZ_t\cdot\vec e_1,~|\vec\mZ_\alpha|\right\rangle\nn\\
&\quad  + 2\pi\gamma\left\langle\vec\mZ\cdot\vec e_1,~(\vec\mZ_t)_\alpha\cdot\vec\mZ_\alpha\,|\vec\mZ_\alpha|^{-1}\right\rangle,\\
\ddt M(\vec\mZ(t)) &= 2\pi\left\langle\vec\mZ\cdot\vec e_1,~\vec\mZ_t\cdot\vec\nu\,|\vec\mZ_\alpha|\right\rangle.
\label{eq:dtM}
\end{align}
\end{subequations}
Now denote 
\begin{equation}
\omega(t) = \frac{\int_{\mR_-(t)}r\,\drz}{\int_{\mR}r\,\drz},\quad t\in[0,T].\label{eq:omega}
\end{equation}
Then it easy to check that $\mX_{_{\mR_-(t)}}-\omega(t)\in\bP$. Choosing $q = \mX_{_{\mR_-(t)}}-\omega(t)$ in \eqref{eq:weak2} and $\zeta = 1$ in \eqref{eq:weak3} and using \eqref{eq:dtM} gives rise to
\begin{align}
\ddt M(\vec\mZ(t)) &=2\pi\left\langle\vec\mZ\cdot\vec e_1,~\vec\mZ_t\cdot\vec\nu\,|\vec\mZ_\alpha|\right\rangle \nn\\
&= 2\pi\left\langle(\vec\mZ\cdot\vec e_1)\vec u\cdot\vec\nu,\,|\vec\mZ_\alpha|\right\rangle = 2\pi\bigl(\nabla_c\cdot[r\,\vec u],~1\bigr)=0.
\label{eq:volumelaw}
\end{align}
Moreover, setting $\vec\chi = \vec u$ in \eqref{eq:weak1}, $q = p$ in \eqref{eq:weak2}, $\zeta= \gamma\varkappa$ in \eqref{eq:weak3} and $\vec\eta=\vec\mZ_t$ in \eqref{eq:weak4}, combining these equations and on recalling \eqref{eq:dtE} yields that
\begin{equation}
\tfrac{1}{2\pi}\ddt E(\rho_{_c}, \vec u(t), \vec\mZ(t)) 
= - 2\Bigl(\mu_{_c}\,r^{-1}\,[\vec u\cdot\vec e_1],~[\vec u\cdot\vec e_1]\Bigr)- 2\Bigl(\mu_{_c}\,r\,\mat{\tD_c}(\vec u),~\mat{\tD_c}(\vec u)\Bigr)+\Bigl(\rho_{_c}r\vec g,~\vec u\Bigr).
\label{eq:energylaw}
\end{equation}
Therefore, we obtain the volume conservation law \eqref{eq:volumelaw} and the energy law \eqref{eq:energylaw} for the weak solution.

On recalling \eqref{eq:kappaEQ}, it is also possible to consider an alternative weak formulation, which treats the curvature of the curve $\Gamma(t)$ as an unknown. To this end, let the initial velocity $\vec u_0\in\bU$ and the initial generating curve $\vec\mZ(\cdot,0)\in\bX$ be given. Then we find $\vec u\in\mathbb{V}$, $p\in L^2(0,T;\bP)$, $\vec\mZ\in H^1(0,T;\bX)$ and $\kappa\in  L^2(0,T;L^2(\bI))$ such that for $t\in(0,T]$
\begin{subequations}\label{eqn:Aweakform}
\begin{align}
&\tfrac{1}{2}\left[\ddt\bigl(\rho_{_c}\vec u,~\vec\chi\,r\bigr) + \bigl(\rho_{_c}\partial_t\vec u,~\vec\chi\,r\bigr) - \bigl(\rho_{_c}\vec u,~\partial_t\vec\chi\,r\bigr)\right]+\mathscr{A}(\rho_{_c},\vec u; \vec u, \vec\chi) \nn\\
&\hspace{1.0cm} + 2\bigl(\mu_{_c}\,r^{-1}\,[\vec u\cdot\vec e_1],~[\vec\chi\cdot\vec e_1]\bigr)+ 2\bigl(\mu_{_c}\,r\,\mat{\tD_c}(\vec u),~\mat{\tD_c}(\vec\chi)\bigr)- \bigl(p,~\bG\cdot[r\,\vec\chi]\bigr)\nn\\
&\hspace{1.0cm} -\,\gamma\big\langle\,(\vec\mZ\cdot\vec e_1)\,\kappa - \vec\nu\cdot\vec e_1,~\vec\nu\cdot\vec\chi\,|\vec\mZ_\alpha|\,\big\rangle = \bigl(\rho_{_c}r\,\vec g,~\vec\chi\bigr)\qquad\forall\vec\chi\in\mathbb{V},\label{eq:Aweak1}\\[0.5em]
\label{eq:Aweak2}
&\hspace{0.2cm}\bigl(\bG\cdot[r\,\vec u], ~q\bigr)=0\qquad\forall q\in \bP,\\[0.5em]
\label{eq:Aweak3}
&\hspace{0.2cm}\big\langle(\vec\mZ\cdot\vec e_1)\,[\partial_t \vec\mZ- \vec u],~\zeta\,\vec\nu\,|\vec\mZ_\alpha|\big\rangle = 0\qquad\forall \zeta\in L^2(\bI),\\[0.5em]
\label{eq:Aweak4}
&\hspace{0.2cm}\big\langle\kappa\,\vec\nu,~\vec\eta\,|\vec\mZ_\alpha|\big\rangle + \big\langle\vec\mZ_\alpha,~\vec\eta_\alpha\,|\vec\mZ_\alpha|^{-1}\big\rangle=0\qquad\forall\vec\eta\in\bX.
\end{align}
\end{subequations}
Similarly, Choosing $q = \mX_{_{\mR_-(t)}}-\omega(t)$ in \eqref{eq:Aweak2} and $\zeta = 1$ in \eqref{eq:Aweak3} then yields that \eqref{eq:volumelaw} is satisfied. Besides, on choosing $\vec\chi=\vec u$ in \eqref{eq:Aweak1}, $q = p$ in \eqref{eq:Aweak2}, $\zeta=\gamma(\kappa - \frac{\vec\nu\cdot\vec e_1}{\vec\mZ\cdot\vec e_1})$ in \eqref{eq:Aweak3} and $\vec\eta = (\vec\mZ\cdot\vec e_1)\vec\mZ_t$ in \eqref{eq:Aweak4}, we obtain \eqref{eq:energylaw} as well.
The advantage of the weak formulation \eqref{eqn:Aweakform} is that it is very
close to \cite[(3.9)]{BGN15stable}, i.e., to the two-dimensional analogue
of \eqref{eqn:3dwf}. Moreover, upon discretization this formulation will admit
an equidistribution property, thanks to \eqref{eq:Aweak4}. However, we will
only be able to prove a stability estimate for discretizations based on the
weak formulation \eqref{eqn:weakform}, utilizing techniques developed in
\cite{BGN19asy,BGN19variational}.

\setcounter{equation}{0}
\section{Unfitted finite element approximations}\label{sec:uffem}

We approximate the weak formulations \eqref{eqn:weakform} and \eqref{eqn:Aweakform} for velocity $\vec u$ and pressure $p$ on the domain $\mR$, 
and for the parameterization $\vec{\mZ}$ and the interfaces mean curvature 
$\varkappa$, or for the curve's curvature $\kappa$, on the reference domain 
$\mathbb{I}$ using the finite element method. 
Let $[0,T]=\bigcup_{m=0}^M[t_m, t_{m+1}]$ be a uniform partition of the time domain with $\ttau = \frac{T}{M}$ and $t_m = m\ttau$ for $m=0,1,\ldots, M$.
 
The reference domain $\bI$ is discretized uniformly into 
$\mathbb{I}=\cup_{j=1}^{J_{\Gamma}}\mathbb{I}_j$ 
with $\mathbb{I}_j=[\alpha_{j-1}, \alpha_j]$, $\alpha_j = jh$ and $h=1/J_{\Gamma}$. We then introduce the finite element spaces
\[
V^h := \Bigl\{\zeta \in C(\overline {\bI}) : \zeta\!\mid_{\bI_j} \nn\
\text{is affine}\ \forall\ j=1,\ldots, J_\Gamma\Bigr\},\qquad
\bX^h:= [V^h]^2\cap\bX.
\]
Let $\{\vec X^m\}_{0\leq m\leq M}$ be an approximation to $\{\vec\mZ(t)\}_{t\in[0,T]}$ with $\vec X^m\in\bX^h$. We  define
$\Gamma^m = \vec X^m\left(\overline{\bI}\right)$, and assume
that
\begin{equation*} 
\vec X^m \cdot\vec e_1 > 0 \quad \text{in }\
\overline {\bI}\setminus \partial \bI
\quad\text{and}\quad
|\vec{X}^m_\alpha| > 0 \quad \text{in }\quad \bI,
\qquad 0\leq m\leq M,
\end{equation*}
so that we have the discrete unit tangential and normal vectors 
\begin{equation*} 
\vec\tau^m = \vec X^m_s = \frac{\vec X^m_\alpha}{|\vec X^m_\alpha|} 
\qquad \mbox{and} \qquad \vec\nu^m = -(\vec\tau^m)^\perp\,.
\end{equation*}

In addition, let $\overline{\mR}=\bigcup_{o\in\mathscr{T}^m} \overline{o}$ 
be a regular partition, where $\mathscr{T}^m:=\bigl\{o_j^m:\,j=1,\ldots,J_\mR^m\bigr\}$ are mutually disjoint and non-degenerate triangles. We introduce the finite element spaces associated with the mesh $\mathscr{T}^m$ as
\begin{align}
S_k^m&:=\left\{\varphi\in C(\overline{\mR}):\; 
\varphi|_{o}\in \mathcal{P}_k(o),\;\forall\ o \in \mathscr{T}^m\right\},\quad k\in\bN_{+},\nn\\
S_0^m&:=\left\{\varphi\in L^2(\mR):\; 
\varphi|_{o}\;\mbox{is constant},\;\forall\ o \in \mathscr{T}^m\right\},\nn
\end{align}
where $\mathcal{P}_k(o)$ denotes the space of polynomials of degree $k$ on $o$. 
Let $\bU^m$ and $\bP^m$ denote the finite element spaces for the numerical solutions of the velocity and pressure, respectively. It is natural to consider the Taylor--Hood pair elements,
\begin{equation}
\label{eq:spaceUP}
\bigl(\mathbb{U}^m,~\mathbb{P}^m\bigr)
=\Bigl([S_2^m]^2\cap\mathbb{U},~
S_1^m \cap\bP\Bigr),
\end{equation}
which guarantees the following inf-sup stability condition for the axisymmetric problem  \cite{Bernardi99,Belhachmi06weighted}: 
\begin{equation}
\inf_{q\in \mathbb{P}^m} 
\sup_{\vec 0\neq\vec\chi\in \mathbb{U}^m}
\frac{\left(q, \bG\cdot[r\,\vec\chi]\right)}
{\norm{r^{\frac{1}{2}}q}\,\norm{\vec\chi}_*}\geq c>0.
\label{eqn:LBB}
\end{equation}
Here $\norm{\cdot}$ denote the usual $L^2$-norm 
on $\mR$, $c$ is a constant and $\norm{\cdot}_*$ is defined as
\begin{equation}
\norm{\vec\chi}_*=\left(\norm{r^{\frac{1}{2}}\bG\chi^r}^2 +\norm{r^{-\frac{1}{2}}\chi^r}^2+ \norm{r^{\frac{1}{2}}\bG\chi^z}^2+\norm{r^\frac{1}{2}\chi^z}^2\right)^{\frac{1}{2}},\nn
\end{equation}
for $\vec\chi = (\chi^r,\chi^z)^T$.

We consider the unfitted mesh approach, meaning that $\mathscr{T}^m$ may not fit the interface, i.e., the line segments making up $\Gamma^m$ are in general not the boundaries of elements in $\mathscr{T}^m$. At $t=t^m$, the interface $\Gamma^m$ divides $\mR$ into two sub-domains $\mR_-^m$ and $\mR_+^m$ with $\mR_-^m$ being the domain enclosed by $\Gamma^m$ and $\partial_z\mR$. 
We split $\mathscr{T}^m$  
into three disjoint subsets $\mathscr{T}_-^m, \mathscr{T}_+^m$ and $\mathscr{T}_\Gamma^m$:
\begin{equation}
\mathscr{T}_\pm^m =\left\{o\in\mathscr{T}^m: o\subset\mR_\pm^m\right\},\quad \mathscr{T}_\Gamma^m =\left\{o\in\mathscr{T}^m: o\cap\Gamma^m\neq\emptyset\right\}.\nn
\end{equation}
Let $\rho^m$ and $\mu^m$ be numerical approximations of
the density $\rho_{_c}(\cdot, t)$ and the viscosity $\mu_{_c}(\cdot, t)$ at $t=t_m$, respectively. 
We define $\rho^m\in S_0^m$ and $\mu^m\in S_0^m$ such that
\begin{equation}
\label{eqn:discreterhovis}
\rho^m|_{o}:=\left\{
\begin{array}{ll}
\rho_-, & \text{if}\ o\in \mathscr{T}_-^m, \vspace{0.15cm}\\
\rho_+, & \text{if}\ o\in\mathscr{T}_+^m,\vspace{0.15cm}\\
\frac{1}{2}(\rho_-+\rho_+), &\text{if}\ o\in \mathscr{T}_{\Gamma}^m,
\end{array}\right.\qquad
\mu^m|_{o}:=\left\{\begin{array}{ll}
\mu_-, &\text{if}\ o\in \mathscr{T}_-^m,\vspace{0.15cm}\\
\mu_+,&\text{if}\ o\in\mathscr{T}_+^m, \vspace{0.15cm}\\
\frac{1}{2}(\mu_-+\mu_+), & \text{if}\ o\in \mathscr{T}_{\Gamma}^m.
\end{array}\right.\nonumber
\end{equation}

In what follows, we denote by $\vec U^m$, $P^m$, $\varkappa^m$, $\kappa^m$ the numerical approximations of $\vec u(\cdot,t)$, $p(\cdot,t)$, $\varkappa(\cdot,t)$ and $\kappa(\cdot, t)$ at time $t_m$, respectively. Besides, we introduce the standard interpolation operators $I_k^m: C(\overline{\mR})\to [S_k^m]^2$ for $k\geq 1$, and the standard projection operator $I_0^m: L^1(\mR)\to S_0^m$ with $(I_0^m\eta)|_o= \frac{1}{|o|}\int_o\eta\drz$ for $o\in\mathscr{T}^m$. Based on suitable discretizations of the weak formulations \eqref{eqn:weakform} and \eqref{eqn:Aweakform}, we next introduce several numerical methods and explore their properties.

\subsection{A linear approximation}

 We start by considering a linear approximation of \eqref{eqn:Aweakform}. Let $\vec U^0\in\bU^0$ be an approximation of the initial velocity $\vec u_0$. Moreover, the generating curve of the initial interface $\Gamma(0)$ is approximated by a polygonal curve $\Gamma^0:=\vec X^0(\overline{\bI})$ with $\vec X^0\in\bX$. We set $\rho^{-1}=\rho^0$. Then for $m\geq0$, we find $\vec U^{m+1}\in \bU^m$, 
$P^{m+1}\in\bP^m$, $\vec X^{m+1}\in \bX^h$ 
and $\kappa^{m+1}\in V^h$ such that  it holds
\begin{subequations}\label{eqn:Afd}
\begin{align}\label{eq:Afd1}
&\hspace{0.1cm}\tfrac{1}{2}\Big[\Bigl(\frac{\rho^m\vec U^{m+1}-I_0^m\rho^{m-1}
I_2^m\vec U^m}{\ttau}, ~\vec\chi\,r\Bigr)
+\Bigl(I_0^m\rho^{m-1}\frac{\vec U^{m+1}-I_2^m\vec U^m}{\ttau}, ~\vec\chi\,r\Bigr)\Big]\nn\\
&\hspace{0.8cm} +\mathscr{A}(\rho^m, I_2^m\vec U^m; \vec U^{m+1},\vec\chi) +2\bigl(\mu^m\,\bD(\vec U^{m+1}),~\bD(\vec\chi)\,r\bigr)\nn\\
&\hspace{0.8cm}+2\bigl(\mu^m\,r^{-1}[\vec U^{m+1}\cdot\vec e_1],~[\vec\chi\cdot\vec e_1]\bigr)^\diamond- \bigl(P^{m+1},~\bG\cdot[r\vec\chi]\bigr)\nn\\
&\hspace{0.8cm} -\,\gamma\,\big\langle(\vec X^{m}\cdot\vec e_1)\,\kappa^{m+1} - \vec\nu^m\cdot\vec e_1,\,\vec\nu^m\cdot\vec\chi\,|\vec X^m_\alpha|\big\rangle  = \big(\rho^m\vec g,~\vec\chi\,r\big),\\[0.3em]
\label{eq:Afd2}
&\hspace{0.1cm}\big(\bG\cdot[r\vec U^{m+1}],~q\Big)=0,\\[0.3em]
\label{eq:Afd3}
&\hspace{0.1cm}\frac{1}{\ttau}\big\langle(\vec X^m\cdot\vec e_1)\,(\vec X^{m+1}-\vec X^m),~\zeta\,\vec\nu^m\,|\vec X^m_\alpha|\big\rangle-\big\langle(\vec X^m\cdot\vec e_1)\,\vec U^{m+1},~\zeta\,\vec\nu^m\,|\vec X^m_\alpha|\big\rangle = 0,\\[0.3em]
\label{eq:Afd4}
&\hspace{0.1cm}\big\langle\kappa^{m+1}\,\vec\nu^m,~\vec\eta\,|\vec X^m_\alpha|\big\rangle^h+ \big\langle\,\vec X^{m+1}_\alpha,~\vec\eta_\alpha\,|\vec X^m_\alpha|^{-1}\big\rangle =0,
\end{align}
\end{subequations}
for all $\bigl(\vec\chi,q, \zeta, \vec\eta\bigr)\in\bU^m\times\bP^m\times V^h\times \bX^h$, where $(\cdot,\cdot)^\diamond$ represents an approximation of the inner product $(\cdot,\cdot)$ that ensures that the corresponding term in \eqref{eq:Afd1} is well defined for all $\vec U^{m+1} \in \bU^m$ and $\vec\chi\in\bU^m$. This is guaranteed by employing quadrature rules with only interior sampling points, such as Gauss--Lobatto quadrature rules. We also introduce $\langle\cdot,\cdot\rangle^h$ as the  mass lumped $L^2$--inner product over $\bI$. In particular,  
for two piecewise continuous functions $\vec v, \vec w$, 
with possible jumps at the nodes $\{\alpha_j\}_{j=1}^{J_\Gamma}$, it is defined via
\[
\big\langle \vec v, \vec u \big\rangle^h = \tfrac{1}{2}\,h\sum_{j=1}^{J_\Gamma} 
\left[(\vec v\cdot\vec u)(\alpha_j^-) + (\vec v \cdot \vec u)(\alpha_{j-1}^+)\right]\quad\mbox{with}\quad g(\alpha_j^\pm)=\underset{\delta\searrow 0}{\lim}\ g(\alpha_j\pm\delta).
\]
The approximation in \eqref{eq:Afd4} leads to an equidistribution property. For a semidiscrete variant this can be made rigorous, while for the fully discrete approximation \eqref{eqn:Afd} the mesh points on the polygonal curve $\Gamma^m$ diffuse tangentially towards an equidistant distribution \cite{BGN19asy}. For a detailed discussion of this equidistribution property, readers are referred to \cite{BGN07,Barrett20}.

\subsection{An unconditionally stable approximation}\label{sec:stabapp}

We next introduce a suitable discretization of \eqref{eqn:weakform}, which leads to an unconditionally stable finite element method. 
With the discrete initial data as before, for $m\geq0$, we find $\vec U^{m+1}\in \bU^m$, 
$P^{m+1}\in\bP^m$, $\vec X^{m+1}\in \bX^h$ 
and $\varkappa^{m+1}\in V^h$ such that it holds
\begin{subequations}\label{eqn:stabfd}
\begin{align}\label{eq:stabfd1}
&\tfrac{1}{2}\Big[\Bigl(\frac{\rho^m\vec U^{m+1}-I_0^m\rho^{m-1}
I_2^m\vec U^m}{\ttau}, ~\vec\chi\,r\Bigr)
+\Bigl(I_0^m\rho^{m-1}\frac{\vec U^{m+1}-I_2^m\vec U^m}{\ttau}, ~\vec\chi\,r\Bigr)\Big]\nn\\
&\hspace{0.8cm} +\mathscr{A}(\rho^m, I_2^m\vec U^m; \vec U^{m+1},\vec\chi)+2\bigl(\mu^m\,\bD(\vec U^{m+1}),~\bD(\vec\chi)\,r\bigr)\nn\\
&\hspace{0.8cm} +2\bigl(\mu^m\,r^{-1}[\vec U^{m+1}\cdot\vec e_1],~[\vec\chi\cdot\vec e_1]\bigr)^\diamond- \,\bigl(P^{m+1},~\bG\cdot[r\vec\chi]\bigr)\nn\\
&\hspace{0.8cm} -\,\gamma\,\big\langle(\vec X^{m}\cdot\vec e_1)\,\varkappa^{m+1}\,\vec\nu^m,\,\vec\chi\,|\vec X^m_\alpha|\big\rangle  = \big(\rho^m\vec g,~\vec\chi\,r\big),\\[0.3em]
\label{eq:stabfd2}
&\big(\bG\cdot[r\vec U^{m+1}],~q\Big)=0,\\[0.3em]
\label{eq:stabfd3}
&\frac{1}{\ttau}\big\langle(\vec X^m\cdot\vec e_1)\,(\vec X^{m+1}-\vec X^m),~\zeta\,\vec\nu^m\,|\vec X^m_\alpha|\big\rangle-\big\langle(\vec X^m\cdot\vec e_1)\,\vec U^{m+1},~\zeta\,\vec\nu^m\,|\vec X^m_\alpha|\big\rangle = 0,\\[0.3em]
\label{eq:stabfd4}
&\big\langle(\vec X^m\cdot\vec e_1)\varkappa^{m+1}\,\vec\nu^m,~\vec\eta\,|\vec X^m_\alpha|\big\rangle + \big\langle\vec\eta\cdot\vec e_1,~|\vec X^{m+1}_\alpha|\big\rangle
+ \big\langle(\vec X^m\cdot\vec e_1)\,\vec X^{m+1}_\alpha,~\vec\eta_\alpha\,|\vec X^m_\alpha|^{-1}\big\rangle =0,
\end{align}
\end{subequations}
for all $\bigl(\vec\chi,q, \zeta, \vec\eta\bigr)\in\bU^m\times\bP^m\times V^h\times \bX^h$. Observe that the only differences with respect to \eqref{eqn:Afd} are
\eqref{eq:stabfd4} and the final term on the left hand side of
\eqref{eq:stabfd1}. 
In addition, we note that \eqref{eqn:stabfd} leads to a system of nonlinear equations due to the presence of $|\vec X^{m+1}_\alpha|$ in the second term of \eqref{eq:stabfd4}. This particular nonlinear approximation contributes to the stability of the interface energy, see \cite{BGN19asy}. We have the following theorem for the introduced method \eqref{eqn:stabfd}, which mimics \eqref{eq:energylaw} on the discrete level.

\begin{theorem}\label{thm:ES} Let $(\vec U^{m+1}, P^{m+1}, \vec X^{m+1}, \varkappa^{m+1})$ be a solution of \eqref{eqn:stabfd}. Then it holds that 
\begin{align}
&\tfrac{1}{2\pi}E(\rho^m,\vec U^{m+1}, \vec X^{m+1}) + 2\ttau\Bigl(\norm{\sqrt{\mu^m\,r}\,\bD(\vec U^{m+1})}^2+\norm{\sqrt{\mu^m\,r^{-1}}[\vec U^{m+1}\cdot\vec e_1]}_\diamond^2\Bigr)\nn\\
&\hspace{1cm}\leq\tfrac{1}{2\pi}E(I_0^m\rho^{m-1},I_2^m\vec U^m,\vec X^m) +\ttau \bigl(\rho^m\,r\vec g,~\vec U^{m+1}\bigr),\quad m = 0,\ldots, M-1,
\label{eq:localES}
\end{align}
where $\norm{\cdot}$ and $\norm{\cdot}_\diamond$ are the induced norms of the inner products $(\cdot,\cdot)$ and $(\cdot,\cdot)^\diamond$, respectively. Moreover, on assuming that 
\begin{equation}
E(I_0^m\rho^{m-1},I_2^m\vec U^m,\vec X^m)\leq E(\rho^{m-1}, \vec U^m, \vec X^m),\quad\mbox{for}\quad m = 0,\ldots, M-1,\label{eq:asp1}
\end{equation}
it holds that for $k=0,\ldots, M-1$.
\begin{align}
&\tfrac{1}{2\pi}E(\rho^k,\vec U^{k+1}, \vec X^{k+1}) + 2\ttau\sum_{m=0}^k\Bigl\{\norm{\sqrt{\mu^m\,r}\,\bD(\vec U^{m+1})}^2+\norm{\sqrt{\mu^m\,r^{-1}}[\vec U^{m+1}\cdot\vec e_1]}_\diamond^2\Bigr\}\nn\\
&\hspace{1cm}\leq\tfrac{1}{2\pi}E(\rho^{0},\vec U^0,\vec X^0) +\ttau \sum_{m=0}^k\bigl(\rho^m\,r\vec g,~\vec U^{m+1}\bigr).\label{eq:globalES}
\end{align}
\end{theorem}
\begin{proof} Setting $\vec\chi=\ttau\,\vec U^{m+1}$ in \eqref{eq:stabfd1}, $q = P^{m+1}$ in \eqref{eq:stabfd2}, $\zeta=\ttau\,\gamma\,\varkappa^{m+1}$ in \eqref{eq:stabfd3} and $\vec\eta = \vec X^{m+1}-\vec X^m$ in \eqref{eq:stabfd4}, and combining these equations yields that
\begin{align}
&\tfrac{1}{2}\left(\rho^m\vec U^{m+1}-I_0^m\rho^{m-1}I_2^m\vec U^m + I_0^m\rho^{m-1}[\vec U^{m+1}-I_2^m\vec U^m], ~\vec U^{m+1}\,r\right) \nn\\
&\hspace{0.8cm}+ 2\ttau\bigl(\mu^m\,r\,\bD(\vec U^{m+1}),~\bD(\vec U^{m+1})\bigr)+2\ttau\bigl(\mu^m\,r^{-1}[\vec U^{m+1}\cdot\vec e_1],~[\vec U^{m+1}\cdot\vec e_1]\bigr)^\diamond\nn\\&\hspace{0.8cm}+\gamma \big\langle(\vec X^{m+1}-\vec X^m)\cdot\vec e_1,~|\vec X^{m+1}_\alpha|\big\rangle+ \gamma\big\langle(\vec X^m\cdot\vec e_1)\,\vec X^{m+1}_\alpha,~(\vec X^{m+1}-\vec X^m)_\alpha\,|\vec X^m_\alpha|^{-1}\big\rangle\nn\\
&\hspace{0.2cm} = \bigl(\rho^m\,r\,\vec g, ~\vec U^{m+1}\bigr).\label{eq:stab1}
\end{align}
It is straightforward to show that the following equation holds
\begin{align}
&\bigl(\rho^m\vec U^{m+1}- I_0^m\rho^{m-1}
I_2^m\vec U^m+ I_0^m\rho^{m-1}[\vec U^{m+1}-I_2^m\vec U^m], ~\vec U^{m+1}\,r\bigr)\nn\\
&\hspace{0.5cm}=\bigl(\rho^m\vec U^{m+1},\vec U^{m+1}r\bigr) - \bigl(I_0^m\rho^{m-1}I_2^m\vec U^m,I_2^m\vec U^m\,r\bigr)\nn\\
&\hspace{2cm}+\bigl(I_0^m\rho^{m-1}[\vec U^{m+1}-I_2^m\vec U^m],[\vec U^{m+1}-I_2^m\vec U^m]\,r\bigr)\nn\\
&\hspace{0.5cm}\geq \bigl(\rho^m\vec U^{m+1},\vec U^{m+1}r\bigr) - \bigl(I_0^m\rho^{m-1}I_2^m\vec U^m,I_2^m\vec U^m\,r\bigr).\label{eq:stab2}
\end{align}
Moreover, using the inequality $\vec a\cdot(\vec a - \vec b)\geq |\vec b|(|\vec a| - |\vec b|)$ for $\vec a, \vec b \in \bR^2$, we have
\begin{align}
&\left\langle (\vec X^{m+1}-\vec X^m) \cdot\vec e_1, ~|\vec X^{m+1}_\alpha|\right\rangle
+ \left\langle (\vec X^m\cdot\vec e_1)\,
\vec X^{m+1}_\alpha,~(\vec X^{m+1} - \vec X^m)_\alpha\, |\vec X^m_\alpha|^{-1} \right\rangle\nn\\
&\hspace{2cm}\geq \left\langle(\vec X^{m+1}-\vec X^m) \cdot\vec e_1, ~|\vec X^{m+1}_\alpha|\right\rangle + \left\langle\vec X^m\cdot\vec e_1 ,~|\vec X^{m+1}_\alpha| - |\vec X^m_{\alpha}| \right\rangle \nn\\
&\hspace{2cm}= \bigl\langle\vec X^{m+1}\cdot\vec e_1,~|\vec X^{m+1}_{\alpha}|\bigr\rangle - \bigl\langle\vec X^m\cdot\vec e_1,~|\vec X^m_\alpha|\bigr\rangle.
\label{eq:stab3}
\end{align}
 Inserting \eqref{eq:stab2} and \eqref{eq:stab3} into \eqref{eq:stab1}, and using \eqref{eq:asyE}, we obtain \eqref{eq:localES} as claimed. Summing \eqref{eq:localES} for $m=0,\ldots, k$ and recalling \eqref{eq:asp1} yields \eqref{eq:globalES}.
\end{proof}

\subsection{Volume-preserving approximations}\label{sec:vpapp}

We recall from \S\ref{sec:25df} that on the continuous level
the volume preservation property \eqref{eq:volumelaw} holds for the two axisymmetric weak formulations \eqref{eqn:weakform} and \eqref{eqn:Aweakform}. 
However, the two schemes \eqref{eqn:Afd} and \eqref{eqn:stabfd} will in general
not conserve the volume of the two phases. In order to motivate a remedy, we
first consider the semidiscrete setting, following similar considerations in
\cite{BGN2013eliminating,BGN15stable}.
In particular, the semidiscrete variant of \eqref{eq:weak2} and 
\eqref{eq:weak3} reads
\begin{subequations}
\begin{align}\label{eq:semi1}
&\hspace{0.2cm}\bigl(\bG\cdot[r\,\vec u^h], ~q^h\bigr)=0\qquad\forall q^h\in \bP^h,\\[0.5em]
&\hspace{0.2cm}\big\langle(\vec\mZ^h\cdot\vec e_1)\,[\partial_t\vec\mZ^h - \vec u^h],~\zeta^h\,\vec\nu^h\,|\vec\mZ_\alpha|\big\rangle = 0\qquad\forall \zeta\in V^h,\label{eq:semi2}
\end{align}
\end{subequations}
where we use the superscript $h$ to denote the corresponding quantities on the semidiscrete level. In a similar manner to \eqref{eq:volumelaw}, and on choosing suitable test functions, it is not difficult to prove that
 \begin{equation}
 \ddt M(\mZ^h(t))=2\pi\big\langle\vec\mZ^h\cdot\vec e_1,~\vec\mZ^h_t\cdot\vec\nu^h\,|\vec\mZ_\alpha|\big\rangle=0
 \end{equation}
if the following condition holds
\begin{equation} \label{eq:extra}
\mX_{\mR_-^h(t)}-\omega^h(t)\in\bP^h\quad\mbox{with}\quad \omega^h(t)= \frac{\int_{\mR_-^h(t)}r\,\drz}{\int_{\mR}r\,\drz},\quad t\in[0,T].
\end{equation}
This suggest to enrich the discrete pressure spaces $\bP^h$ with the single
extra function in \eqref{eq:extra}. On the fully discrete level this
corresponds to ensuring that
\begin{equation}\label{eq:enrichp}
\mX_{_{\mR_-^m}}- \omega^m\in\bP^m\quad{\rm with} \quad \omega^m =\frac{\int_{\mR_-^m}r\,\drz}{\int_\mR r\,\drz},\quad m = 0,\ldots, M-1.
\end{equation}
Therefore we enforce this by adding the single basis function $\mX_{_{\mR_-^m}}-\omega^m$
to the discrete pressure space $\bP^m$, a procedure which we call XFEM. 
The contribution of this new basis function to \eqref{eq:Afd1}, \eqref{eq:Afd2}, \eqref{eq:stabfd1} and \eqref{eq:stabfd2} can be rewritten in terms of integrals over $\bI$, on noting that
\begin{equation}
\bigl(\bG\cdot[r\,\vec\chi],~\mX_{_{\mR_-^m}}\bigr) 
= \int_{\mR_-^m}\bG\cdot[r\,\vec\chi]\,\drz 
= \bigl\langle(\vec X^m\cdot\vec e_1)\,\vec\chi,~\vec\nu^m\bigr\rangle\qquad \forall\vec\chi\in\bU^m.
\end{equation}
We note that this XFEM technique was first introduced for the discretizations of \eqref{eqn:3dwf} in \cite{BGN15stable}, leading to good volume conservation in practice. Nevertheless, in order to satisfy an exact volume preservation also on the fully discrete level, it turns out that suitable time-weighted approximations of the interface normals are needed, see \cite{BZ21SPFEM,BGNZ22volume}.

 We follow the idea in \cite{BGNZ22volume} and introduce $\vec f^{m+\frac{1}{2}}$ as an appropriate approximation of the quantity $\vec f=\vec\mZ\cdot\vec e_1\,|\vec \mZ_\alpha|\,\vec\nu$. Precisely, we define $\vec f^{m+\frac{1}{2}} \in [L^\infty(\bI)]^2$ by
\begin{align}
\label{eq:averagenor}
\vec f^{m+\frac{1}{2}}= -\tfrac{1}{6}\bigl[(\vec X^m\cdot\vec e_1)\,\vec X^m_{\alpha}+4(\vec X^{m+\frac{1}{2}}\cdot\vec e_1)\vec X^{m+\frac{1}{2}}_\alpha+(\vec X^{m+1}\cdot\vec e_1)\,\vec X^{m+1}_{\alpha}\bigr]^\perp,
\end{align}
where $\vec X^{m+\frac{1}{2}}=\frac{1}{2}(\vec X^m + \vec X^{m+1})$. Then we have the following result from \cite[Lemma 3.1]{BGNZ22volume}. For the sake of completeness, we include the proof here as well.

\begin{lemma}\label{lem:vol}
Let $\vec X^m \in\bX^h$ and $\vec X^{m+1} \in \bX^h$. Then it holds that
\begin{equation} \label{eq:MMf}
M(\vec X^{m+1})-M(\vec X^m) = 2\pi\,\bigl\langle\vec X^{m+1}-\vec X^m,~\vec f^{m+\frac{1}{2}}\bigr\rangle.
\end{equation}
\end{lemma} 
\begin{proof}
We first introduce a linear approximation in time between $\vec X^m$ and $\vec X^{m+1}$ as
\begin{equation} \label{eq:Xht}
\vec X^h(t):=\frac{t_{m+1} - t}{\ttau}\vec X^{m+1} + \frac{t-t_m}{\ttau}\vec X^m\in \bX^h,\quad t\in[t_m, t_{m+1}].
\end{equation}
Then applying \eqref{eq:dtM} to $\vec X^h(t)$ yields that
\begin{equation}
\ddt M(\vec X^h(t))=2\pi\left\langle\vec X^h\cdot\vec e_1,~\vec X^h_t\cdot\vec\nu^h\,|\vec X^h_\alpha|\right\rangle \quad \forall t\in(t_m, t_{m+1}).\label{eq:dtmh}
\end{equation}
Integrating \eqref{eq:dtmh} from $t_m$ to $t_{m+1}$, and noting that 
$\vec X^h_t = \frac{\vec X^{m+1}-\vec X^m}{\ttau}$ and 
$\vec\nu^h = -\frac{(\vec X^h_\alpha)^\perp}{|\vec X^h_\alpha|}$,
yields that
\begin{align*}
M(\vec X^{m+1})-M(\vec X^m) &= 2\pi\int_{t_m}^{t_{m+1}}\int_\bI(\vec X^h\cdot\vec e_1)\,\frac{\vec X^{m+1}-\vec X^m}{\ttau}\cdot(-\vec X^h_\alpha)^\perp\,\rd\alpha\rd t\nn\\
&=2\pi\int_\bI\frac{\vec X^{m+1}-\vec X^m}{\ttau}\cdot\int_{t_m}^{t_{m+1}}(\vec X^h\cdot\vec e_1)(-\vec X^h_\alpha)^\perp\rd t\rd\alpha\nn\\
&=2\pi\int_\bI (\vec X^{m+1}-\vec X^m)\cdot\vec f^{m+\frac{1}{2}}\rd\alpha,
\end{align*}
where in the last equality we employed Simpson's rule on noting that the integrand $(\vec X^h\cdot\vec e_1)\,(-\vec X^h_\alpha)^\perp$ is a quadratic function in $t$.
\end{proof}

We are now in a position to adapt the numerical method \eqref{eqn:stabfd} to yield an exactly volume conserving scheme as follows. For $m\geq0$, we find $\vec U^{m+1}\in \bU^m$, $P^{m+1}\in\bP^m$, $\vec X^{m+1}\in \bX^h$ 
and $\varkappa^{m+1}\in V^h$ such that it holds 
\begin{subequations}\label{eqn:spfd}
\begin{align}\label{eq:spfd1}
&\tfrac{1}{2}\Big[\Bigl(\frac{\rho^m\vec U^{m+1}-I_0^m\rho^{m-1}
I_2^m\vec U^m}{\ttau}, ~\vec\chi\,r\Bigr)
+\Bigl(I_0^m\rho^{m-1}\frac{\vec U^{m+1}-I_2^m\vec U^m}{\ttau}, ~\vec\chi\,r\Bigr)\Big]\nn\\
&\hspace{0.8cm} +\mathscr{A}(\rho^m, I_2^m\vec U^m; \vec U^{m+1},\vec\chi) +2\bigl(\mu^m\,\bD(\vec U^{m+1}),~\bD(\vec\chi)\,r\bigr)\nn\\
&\hspace{0.8cm}+2\bigl(\mu^m\,r^{-1}[\vec U^{m+1}\cdot\vec e_1],~[\vec\chi\cdot\vec e_1]\bigr)^\diamond- \,\bigl(P^{m+1},~\bG\cdot[r\vec\chi]\bigr)\nn\\
&\hspace{0.8cm} -\,\gamma\,\big\langle(\vec X^{m}\cdot\vec e_1)\,\varkappa^{m+1},\,\vec\nu^m\cdot\vec\chi\,|\vec X^m_\alpha|\big\rangle  = \big(\rho^m\vec g,~\vec\chi\,r\big),\\[0.3em]
\label{eq:spfd2}
&\hspace{0.2cm}\big(\bG\cdot[r\vec U^{m+1}],~q\Big)=0,\\[0.3em]
\label{eq:spfd3}
&\hspace{0.2cm}\frac{1}{\ttau}\big\langle\vec X^{m+1}-\vec X^m,~\zeta\,\vec f^{m+\frac{1}{2}}\big\rangle-\big\langle(\vec X^m\cdot\vec e_1)\,\vec U^{m+1},~\zeta\,\vec\nu^m\,|\vec X^m_\alpha|\big\rangle = 0,\\[0.3em]
\label{eq:spfd4}
&\hspace{0.2cm}\big\langle\varkappa^{m+1}\,\vec f^{m+\frac{1}{2}},~\vec\eta\big\rangle + \big\langle\vec\eta\cdot\vec e_1,~|\vec X^{m+1}_\alpha|\big\rangle+ \big\langle(\vec X^m\cdot\vec e_1)\,\vec X^{m+1}_\alpha,~\vec\eta_\alpha\,|\vec X^m_\alpha|^{-1}\big\rangle =0,
\end{align}
\end{subequations}
for all $\left(\vec\chi, q, \zeta, \vec\eta\right)\in\bU^m\times\bP^m\times V^h\times \bX^h$.  Here in \eqref{eq:spfd3}, we employ a semi-implicit treatment $\vec f^{m+\frac{1}{2}}$ to approximate $(\vec\mZ\cdot\vec e_1)|\vec\mZ_\alpha|\,\vec\nu$ instead of the explicit treatment in \eqref{eq:stabfd3}. This then enables the exact volume preservation. Besides, we also employ $\vec f^{m+\frac{1}{2}}$ in \eqref{eq:spfd4} in order to maintain our stability bound. We then have the following theorem for the scheme \eqref{eqn:spfd}, which mimics \eqref{eq:volumelaw} and \eqref{eq:energylaw} on the discrete level. 

\begin{theorem}\label{thm:SP}
Let $(\vec U^{m+1}, P^{m+1}, \vec X^{m+1}, \varkappa^{m+1})$ be a solution of \eqref{eqn:spfd}. Then \eqref{eq:localES} holds.
If the condition \eqref{eq:asp1} is satisfied, then 
\eqref{eq:globalES} holds for $k=0,\ldots, M-1$.
Moreover, if \eqref{eq:enrichp} holds then
\begin{equation}\label{eq:vol}
M(\vec X^{m+1})=M(\vec X^m),\quad m = 0,\ldots, M-1.
\end{equation}
\end{theorem}
\begin{proof}
The stability results can be shown as in the proof of Theorem~\ref{thm:ES}.
For the volume preservation, we choose $q = (\mX_{_{\mR_-^m}} - \omega^m)\in \bP^m$ in \eqref{eq:spfd2} and obtain
\begin{equation}
\label{eq:vol1}
\bigl(\bG\cdot[r\,\vec U^{m+1}], \mX_{_{\mR_-^m}}\bigr) - \omega^m\bigl(\bG\cdot[r\,\vec U^{m+1}],1\bigr)=0.
\end{equation}
Using integration by parts, we can recast \eqref{eq:vol1} as
\begin{align}
0&=\bigl(\bG\cdot[r\,\vec U^{m+1}), \mX_{_{\mR_-^m}}\bigr)-\omega^m\bigl(\bG\cdot[r\,\vec U^{m+1}],1\bigr)\nn\\
 &=\int_{\Gamma^m}(\vec X^m\cdot\vec e_1)\,\vec U^{m+1}\cdot\vec \nu^m\,\rd s = \big\langle(\vec X^m\cdot\vec e_1)\,\vec U^{m+1},~\vec\nu^m|\vec X^m_\alpha|\big\rangle.
 \label{eq:vol2}
\end{align}
On the other hand, setting $\zeta=\ttau$ in \eqref{eq:spfd3} and using \eqref{eq:vol2}  we obtain
\begin{equation}
\bigl<\vec X^{m+1}-\vec X^m,~\vec f^{m+\frac{1}{2}}\bigr> =\ttau\big\langle(\vec X^m\cdot\vec e_1)\,\vec U^{m+1},~\vec\nu^m|\vec X^m_\alpha|\big\rangle= 0,
\end{equation}
which implies the desired volume preservation \eqref{eq:vol} on recalling Lemma \ref{lem:vol}.
\end{proof}

It is not difficult to also adapt the linear scheme \eqref{eqn:Afd} 
into an exactly volume preserving method. In fact, the new scheme,
which is now nonlinear, is given by \eqref{eqn:Afd} with \eqref{eq:Afd3} replaced by \eqref{eq:spfd3}. 
\begin{theorem}\label{thm:vc}
Let $(\vec U^{m+1}, P^{m+1},~\vec X^{m+1},~\kappa^{m+1})$ be a solution of the adapted scheme \eqref{eqn:Afd} with \eqref{eq:Afd3} replaced by \eqref{eq:spfd3}. On assuming that \eqref{eq:enrichp} is satisfied, it holds that
\begin{equation}
\label{eq:volnew}
M(\vec X^{m+1}) = M(\vec X^m),\quad m= 0,\ldots, M-1.
\end{equation}
\end{theorem}
\begin{proof}
The proof is exactly the same as that of \eqref{eq:vol}.
\end{proof}

\subsection{Nonlinear solvers} 
We note that apart from \eqref{eqn:Afd}, all our introduced schemes lead to 
systems of nonlinear equations at each time step. 
The linear systems arising from \eqref{eqn:Afd} can be solved with the help 
of a Schur complement approach and preconditioned Krylov iterative solvers,
as described in \cite{BGN15stable}. 
For the nonlinear systems resulting from the other schemes, we employ a 
Picard-type iteration. For example, in order to solve for \eqref{eqn:spfd}, we first set the initial guess $\vec X^{m+1,0}=\vec X^m$. Then, for each $l\geq 0$, we find  $(\vec U^{m+1,l+1}, P^{m+1,l+1}, \vec X^{m+1,l+1}, \varkappa^{m+1, l+1})\in\bU^m\times\bP^m\times\bX^h\times V^h$ such that it holds
\begin{subequations}\label{eqn:picardfd}
\begin{align}\label{eq:picardfd1}
&\tfrac{1}{2}\Big[\Bigl(\frac{\rho^m\vec U^{m+1, l+1}-I_0^m\rho^{m-1}
I_2^m\vec U^m}{\ttau}, ~\vec\chi\,r\Bigr)
+\Bigl(I_0^m\rho^{m-1}\frac{\vec U^{m+1, l+1}-I_2^m\vec U^m}{\ttau}, ~\vec\chi\,r\Bigr)\Big] \nn\\
&\hspace{0.8cm}
+\mathscr{A}(\rho^m, I_2^m\vec U^m; \vec U^{m+1, l+1},\vec\chi) +2\bigl(\mu^m\,\bD(\vec U^{m+1, l+1}),~\bD(\vec\chi)\,r\bigr)\nn\\
&\hspace{0.8cm}+2\bigl(\mu^m\,r^{-1}[\vec U^{m+1, l+1}\cdot\vec e_1],~[\vec\chi\cdot\vec e_1]\bigr)^\diamond- \,\bigl(P^{m+1, l+1},~\bG\cdot[r\vec\chi]\bigr)\nn\\
&\hspace{0.8cm} -\,\gamma\,\big\langle(\vec X^{m}\cdot\vec e_1)\,\varkappa^{m+1, l+1}\,\vec\nu^m,\,\vec\chi\,|\vec X^m_\alpha|\big\rangle  = \big(\rho^m\vec g,~\vec\chi\,r\big),\\[0.3em]
\label{eq:picardfd2}
&\hspace{0.1cm}\big(\bG\cdot[r\vec U^{m+1, l+1}],~q\Big)=0,\\[0.3em]
\label{eq:picardfd3}
&\hspace{0.2cm}\frac{1}{\ttau}\big\langle\vec X^{m+1, l+1}-\vec X^m, ~\zeta\vec f^{m+\frac{1}{2}, l}\big\rangle-\big\langle(\vec X^m\cdot\vec e_1)\,\vec U^{m+1, l+1},~\zeta\,\vec\nu^m\,|\vec X^m_\alpha|\big\rangle = 0,\\[0.3em]
\label{eq:picardfd4}
&\hspace{0.2cm}\big\langle\varkappa^{m+1, l+1}\,\vec f^{m+\frac{1}{2}, l},~\vec\eta\big\rangle + \big\langle\vec\eta\cdot\vec e_1,~|\vec X^{m+1, l}_\alpha|\big\rangle
+ \big\langle(\vec X^m\cdot\vec e_1)\,\vec X^{m+1, l+1}_\alpha,~\vec\eta_\alpha\,|\vec X^m_\alpha|^{-1}\big\rangle =0,
\end{align}
\end{subequations}
for all $\bigl(\vec\chi,q, \zeta, \vec\eta\bigr)\in\bU^m\times\bP^m\times V^h\times \bX^h$, where $\vec f^{m+\frac{1}{2}, l}$ is a lagged approximation which follows \eqref{eq:averagenor} directly except that $\vec X^{m+1}$ was replaced by $\vec X^{m+1, l}$. 
The system \eqref{eqn:picardfd} is linear and of the same general structure as
\eqref{eqn:Afd}. Hence it can be efficiently solved with the aforementioned
solution techniques. We repeat the above iteration until
\begin{equation}
\max_{0\leq j\leq J_\Gamma}\left|\vec X^{m+1,l+1}(\alpha_j) - \vec X^{m+1, l}(\alpha_j)\right|\leq {\rm tol},
\end{equation}
where ${\rm tol}$ is a chosen tolerance.

\setcounter{equation}{0}
\section{Numerical results}\label{sec:num}

We implemented our unfitted finite element approximations within the 
finite element toolbox Alberta, \cite{Alberta}. 
For the quadrature $(\cdot,\cdot)^\diamond$ we employ a Gauss--Lobatto
formula that is exact for polynomials up to degree $17$.
For the velocity and pressure approximation in the bulk we always use
the lowest order Taylor--Hood element, together with our XFEM extension of the
pressure space, recall \eqref{eq:spaceUP} and \S\ref{sec:vpapp}.

As pointed out in \cite[p.~125]{Ganesan08accurate}, it is often quite
straightforward to extend existing 2d finite element codes to the axisymmetric
setting. For example, we use the same adaptive mesh strategies as in
\cite{BGN15stable} that lead to locally refined bulk meshes close to the
interface. For the specification of the mesh parameters we use the same
notation as there: $n\adapt{k}{l}$ means the time step size
is set to $\ttau = 10^{-3}/n$, the discrete interface is partitioned into
$J_\Gamma=2^k$ elements and the two bulk mesh refinement parameters
$N_f=2^k$ and $N_c=2^l$ indicate the mesh sizes employed close to the interface
and far away from it, respectively. We refer to \cite{BGN15stable} for a more
detailed description of the adaptive mesh refinement strategies.
For simplicity, we denote the introduced methods \eqref{eqn:stabfd} and 
\eqref{eqn:spfd} by $\StabN^h$ and $\StabV^h$, respectively. Moreover, 
we use $\EquiL^h$ to denote the method \eqref{eqn:Afd}, and $\EquiV^h$ to 
denote the scheme \eqref{eqn:Afd} with \eqref{eq:Afd3} 
replaced by \eqref{eq:spfd3}.

\subsection{The rising bubble}

Inspired by the 2d benchmark computations proposed in \cite{Hysing2009}, see 
\cite{BGN15stable} for their generalizations to 3d, we study the dynamics of a rising bubble in a bounded cylinder of diameter $1$ and height $2$. The initial interface of the bubble is given by a sphere of radius $\frac{1}{4}$ and centre $(0,0,\frac{1}{2})^T$. This means that in cylindrical coordinates, we have $\mR = [0, \frac{1}{2}]\times[0,2]$ and $\Gamma(0)=\left\{\vec\mZ\in\mR:\;|\vec\mZ-(0, \frac{1}{2})^T|=\frac{1}{4}\right\}$. In order to monitor the progress of the rising bubble, we introduce the discrete benchmark quantities
\begin{subequations}
\begin{alignat}{3}
\cira|_{t_m} &:= \left(\frac{9}{2\pi^2}\right)^{\frac{1}{3}}\frac{[M(\vec X^m)]^{\frac{2}{3}}}{\int_0^1(\vec X^m\cdot\vec e_1)|\vec X^m_\alpha|\rd\alpha},\qquad &&V_c|_{t_m}:=\frac{2\pi\int_{\mR_-^m}(\vec U^m\cdot\vec e_2)r\,\drz}{M(\vec X^m)},\nn\\
&z_c|_{t_m}:=\frac{2\pi\int_{\mR_-^m}(\vec\id\cdot\vec e_2)r\,\drz}{M(\vec X^m)},\qquad 
&&v_\Delta|_{t_m}: = \frac{M(\vec X^m)-M(\vec X^0)}{M(\vec X^0)},\nn
\end{alignat}
\end{subequations}
where $\cira$ denotes the degree of sphericity of the bubble, $V_c$ is the bubble's rise velocity, $z_c$ is the bubble's centre of mass in the vertical direction,  $v_\Delta$ is the relative volume loss.

\begin{table}[!htp]
\centering
\def\temptablewidth{0.85\textwidth}
\vspace{0pt}
\caption{Benchmark quantities of the rising bubble in case I.}\label{tab:case1}
{\rule{\temptablewidth}{1pt}}
\begin{tabular}{c|ccc|ccc}
&\multicolumn{3}{c|}{$\StabN^h$} &\multicolumn{3}{c}{$\StabV^h$}\\ \hline 
  &${\rm adapt}_{5,3}$ & ${\rm adapt}_{7,3}$   &$2{\rm adapt}_{9,4}$  &${\rm adapt}_{5,3}$ & ${\rm adapt}_{7,3}$   &$2{\rm adapt}_{9,4}$ \\\hline  
$\cira_{\min}$ &0.9630  &0.9536 &0.9501  &0.9630 &0.9536   &0.9501 \\\hline 
$t_{\cira = \cira_{\min}}$ &3.000 &2.9960 &3.0000 &3.000  & 3.0000  &3.0000 \\[0.4em]\hline 
$V_{c,\max}$ &0.3687 &0.3637 &0.3642  &0.3686 & 0.3640 &0.3643    \\\hline 
$t_{_{V_c = V_{c,\max}}}$ &1.233 &0.8950 &0.9255 &1.233  & 0.8980  &0.9255   \\\hline 
$z_c(t=3)$&1.4834  & 1.4848 &  1.4897  &1.4835  &1.4848 &1.4897  \\\hline  
 $v_\Delta(t=3)$&8.50E-5 & -8.27E-5 &-2.83E-5   &0 &0 &0 
 \end{tabular}
{\rule{\temptablewidth}{1pt}}
{\rule{\temptablewidth}{1pt}}
\begin{tabular}{c|ccc|ccc}
&\multicolumn{3}{c|}{$\EquiL^h$} &\multicolumn{3}{c}{$\EquiV^h$}\\ \hline 
 &${\rm adapt}_{5,2}$ & ${\rm adapt}_{7,3}$   &$2{\rm adapt}_{9,4}$ &${\rm adapt}_{5,2}$  & ${\rm adapt}_{7,3}$   &$2{\rm adapt}_{9,4}$    \\\hline  
$\cira_{\min}$ &0.9603 & 0.9534 &0.9501  &0.9603 &0.9534   &0.9501  \\\hline 
$t_{\cira = \cira_{\min}}$ &3.0000 &3.0000  &3.0000 & 3.0000  &3.0000 &3.0000\\[0.4em]\hline 
$V_{c,\max}$&0.3681 &0.3637 &0.3642  &0.3683  & 0.3637 &0.3642    \\\hline 
$t_{_{V_c = V_{c,\max}}}$ &1.236  &0.8980 &0.9260  &1.236 & 0.8980  &0.9260  \\\hline 
$z_c(t=3)$ &1.4815 &1.4846 &  1.4896  &1.4815 &1.4846 &1.4896 
\\\hline  
 $v_\Delta(t=3)$ &-4.39E-5 &-4.58E-5 & -2.54E-5  &0 &0 &0 
   \end{tabular}
{\rule{\temptablewidth}{1pt}}
\end{table}

\begin{figure}[!htp]
\centering
\includegraphics[width=0.7\textwidth]{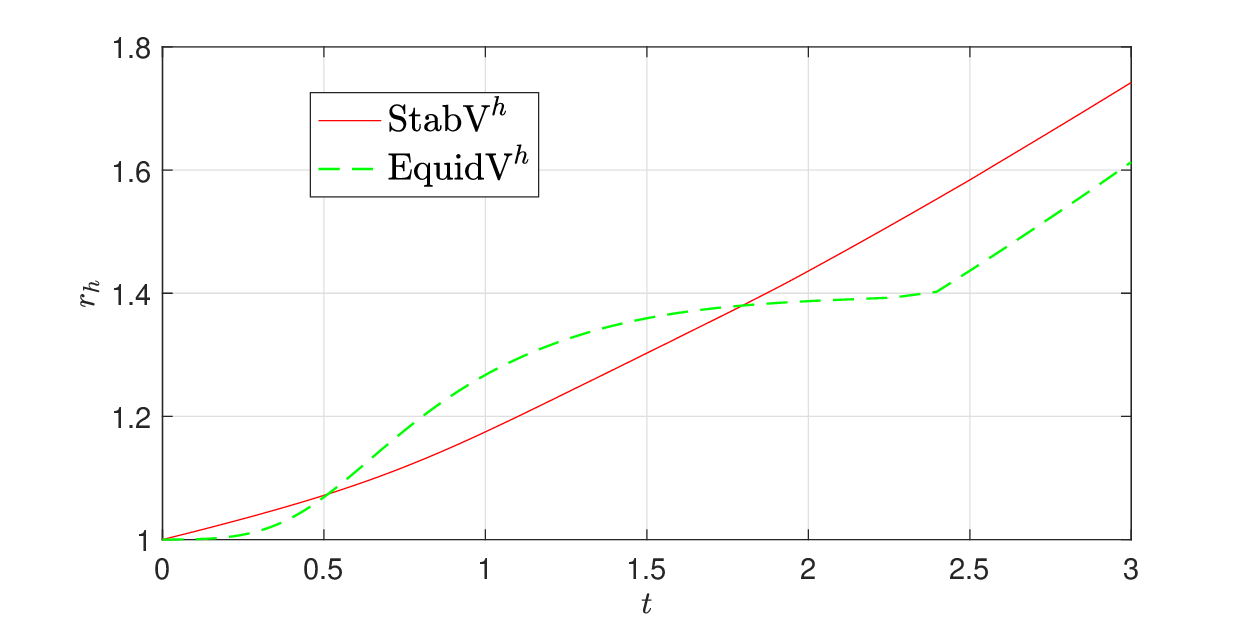}
\caption{ The time history of the mesh ratio function $r_h$ defined in \eqref{eq:rh} for the methods $\StabV^h$ and $\EquiV^h$ with $2{\rm adapt}_{9,4}$. }
\label{fig:meshQ}
\end{figure}

\begin{figure}[t]
\centering
\includegraphics[width=0.45\textwidth]{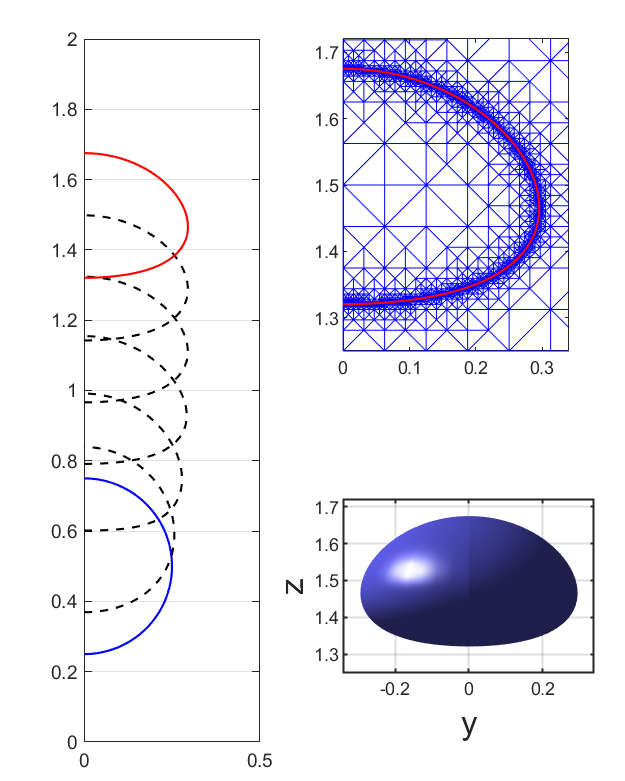}
\includegraphics[width=0.45\textwidth]{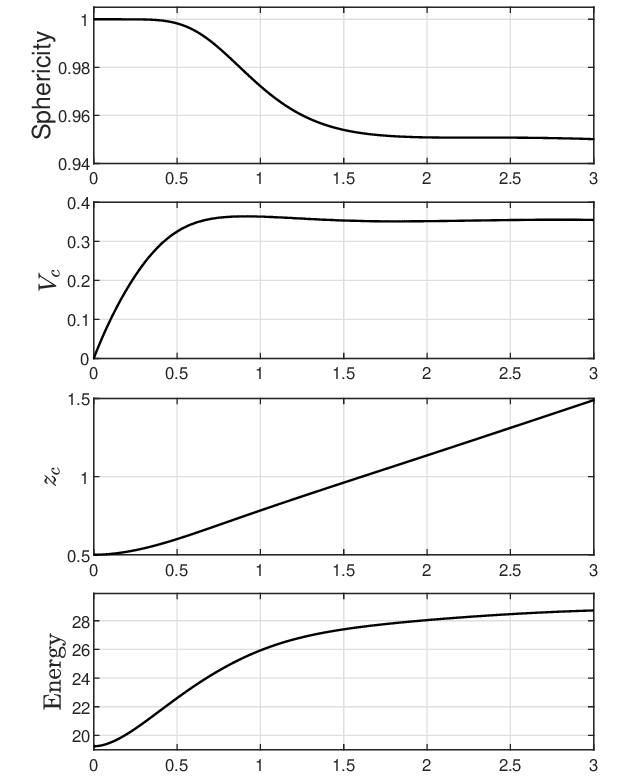}
\caption{Evolution of the rising bubble in case I using the method $\EquiV^h$ with $2{\rm adapt}_{9,4}$. On the left we show the generating curves of the interfaces at $t=0, 0.5, \ldots, 3$ with visualisations of the computational mesh and the axisymmetric surface $\mS(t)$ at the final time $t=3$.  On the right are plots of the discrete benchmark quantities versus time.}
\label{fig:case1}
\end{figure}

We consider two difference cases, and the physical parameters are given as in \cite{Hysing2009}, see also \cite{BGN15stable}.
\begin{itemize}
\item Case I:
\begin{equation}
\rho_+=1000, \quad\rho_- = 100, \quad\mu_+=10,\quad\mu_-=1,\quad \gamma = 24.5,\quad \vec g = (0, -0.98)^T.\nn
\end{equation}
\item Case II:
\begin{equation}
\rho_+=1000, \quad\rho_- = 1, \quad\mu_+=10,\quad\mu_-=0.1,\quad \gamma = 1.96,\quad \vec g = (0, -0.98)^T.\nn
\end{equation}
\end{itemize}
We first consider case I and use the P2-P1 pair elements in \eqref{eq:spaceUP} with XFEM that discussed in \S\ref{sec:vpapp} to improve the volume preservation. The benchmark results computed by the four introduced methods are reported in Table \ref{tab:case1}. We observe that these methods can produce quite similar results, and the volume of the bubble is well preserved. In particular, we observe the exact volume preservation for $\StabV^h$ and $\EquiV^h$, which numerically verifies \eqref{eq:vol} and \eqref{eq:volnew}. To examine the quality of the interface mesh, we introduce the mesh ratio function 
\begin{equation}
r_h|_{t_m} = \frac{\max_{1\leq j\leq J_\Gamma}|\vec X^m(\alpha_j) - \vec X^m(\alpha_{j-1}|}{\min_{1\leq j\leq J_\Gamma}|\vec X^m(\alpha_j) - \vec X^m(\alpha_{j-1}|},\label{eq:rh}
\end{equation}
and plot the time history of $r_h$ in Fig.~\ref{fig:meshQ}. We find that for both the methods $\StabV^h$ and $\EquiV^h$, the value of the mesh ratio function will in general not exceed 2.0 for $t\leq 3$. This implies the good mesh quality. For the experiment using the method $\EquiV^h$ with $2{\rm adapt}_{9,4}$, several snapshots of the generating curves are shown in Fig.~\ref{fig:case1}. Moreover, in Fig.~\ref{fig:case1} we visualise the computational mesh and the bubble at the final time $t=3$. Plots of the benchmark quantities versus time are shown as well in Fig.~\ref{fig:case1}.

\begin{table}[!htp]
\centering
\def\temptablewidth{0.85\textwidth}
\vspace{0pt}
\caption{Benchmark quantities of the rising bubble in case II.}\label{tab:case2}
{\rule{\temptablewidth}{1pt}}
\begin{tabular}{c|ccc|ccc}
&\multicolumn{3}{c|}{$\StabN^h$} &\multicolumn{3}{c}{$\StabV^h$}\\ \hline 
 &${\rm adapt}_{5,2}$  & ${\rm adapt}_{7,3}$   &$2{\rm adapt}_{9,4}$  &${\rm adapt}_{5,2}$  & ${\rm adapt}_{7,3}$   &$2{\rm adapt}_{9,4}$   \\\hline  
$\cira_{\min}$ &0.8117 & 0.7872 &0.7682 &0.8108 &0.7873   &0.7683  \\\hline 
$t_{\cira = \cira_{\min}}$ &1.5000 &1.5000  &1.5000 & 1.5000  &1.5000 &1.5000\\[0.4em]\hline 
$V_{c,\max}$ &0.3621 &0.3693 &0.3717  &0.3633  & 0.3693 &0.3717   \\\hline 
$t_{_{V_c = V_{c,\max}}}$&0.6970 &0.5640 &0.5440  &0.6180 & 0.5650  &0.5440   \\\hline 
$z_c(t=1.5)$ &0.9833 & 0.9901 &0.9972  &0.9831 &0.9900 &0.9972  \\\hline  
 $v_\Delta(t=1.5)$ &7.94E-7 &-2.68E-4 & -1.17E-4   &0 &0 &0 
  \end{tabular}
{\rule{\temptablewidth}{1pt}}
{\rule{\temptablewidth}{1pt}}
\begin{tabular}{c|ccc|ccc}
&\multicolumn{3}{c|}{$\EquiL^h$} &\multicolumn{3}{c}{$\EquiV^h$}\\ \hline 
&${\rm adapt}_{5,2}$  & ${\rm adapt}_{7,3}$   &$2{\rm adapt}_{9,4}$  &${\rm adapt}_{5,2}$  & ${\rm adapt}_{7,3}$   &$2{\rm adapt}_{9,4}$\\\hline  
$\cira_{\min}$ &0.8014 & 0.7865 &0.7682  &0.8013 &0.7866   &0.7682  \\\hline 
$t_{\cira = \cira_{\min}}$ &1.5000 &1.5000  &1.5000 & 1.5000  &1.5000 &1.5000\\[0.4em]\hline 
$V_{c,\max}$ &0.3630 &0.3694 &0.3717  &0.3629 & 0.3693 &0.3717    \\\hline 
$t_{_{V_c = V_{c,\max}}}$ &0.6190 &0.5640 &0.5440  &0.6190 & 0.5640  &0.5440  \\\hline 
$z_c(t=1.5)$ &0.9817 &0.9901 & 0.9972  &0.9817 &0.9901 &0.9972  \\\hline  
 $v_\Delta(t=1.5)$&1.25E-5  & -1.15E-4 & -9.95E-5  &0 &0 &0 
  \end{tabular}
{\rule{\temptablewidth}{1pt}}
\end{table}

\begin{figure}[!htp]
\centering
\includegraphics[width=0.95\textwidth]{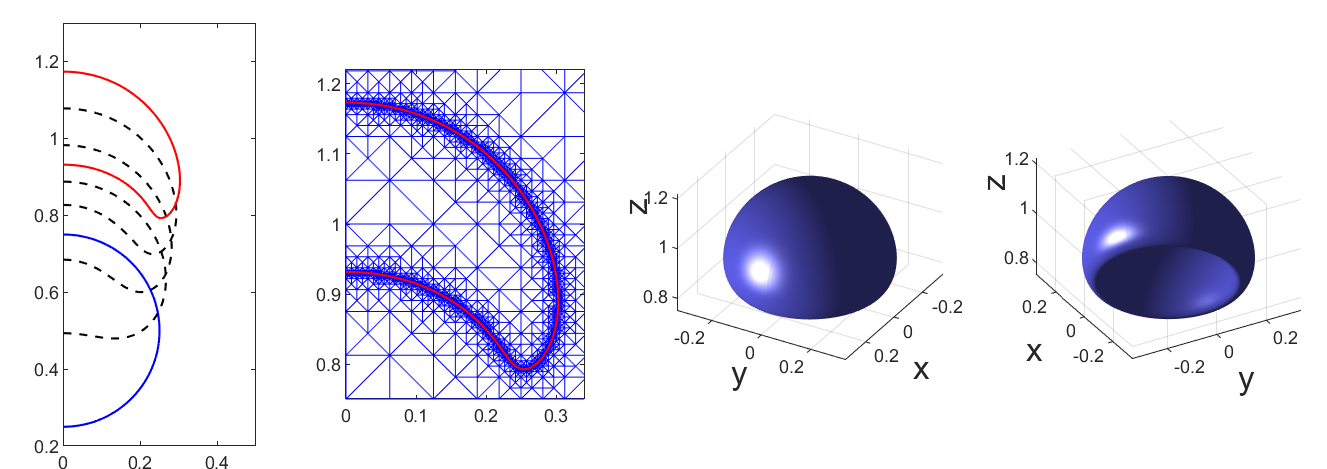}
\caption{Evolution of the rising bubble in case II using the method $\StabV^h$ with $2{\rm adapt}_{9,4}$, where we show the generating curves of the interfaces at $t=0, 0.3, 0.6, 1.2, 1.5$, the visualisations of the computational mesh at the final time,  and the axisymmetric surface $\mS(t)$ at the final time $t=1.5$ with sights in different directions.}
\label{fig:case2}
\end{figure}

\begin{figure}[!htp]
\centering
\includegraphics[width=0.9\textwidth]{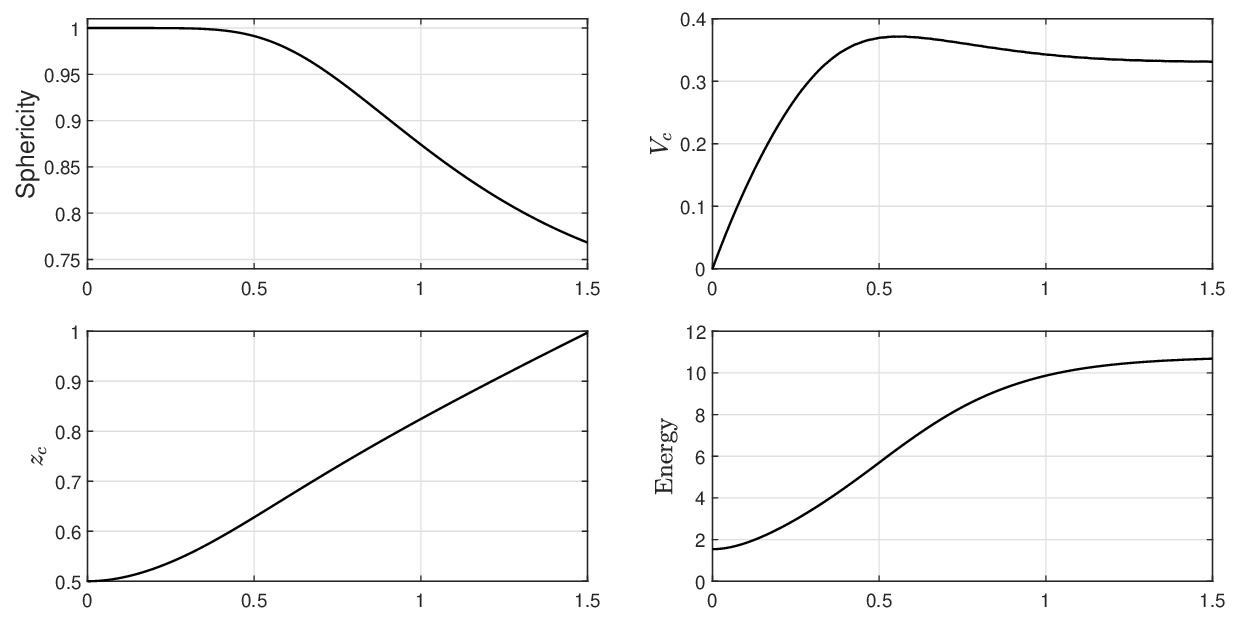}
\caption{The time evolution of the discrete benchmark quantities in Fig.~\ref{fig:case2}.}
\label{fig:case2Q}
\end{figure}

We next consider the rising bubble problem in case II. In this case, the rising bubble can exhibit strong deformations due to the high density ratio and viscosity ratio. The benchmark results are reported in Table \ref{tab:case2}. For the experiment using the method $\StabV^h$ with $2{\rm adapt}_{9,4}$, we show the evolving interfaces and the plots of the benchmark results in Fig.~\ref{fig:case2} and Fig.~\ref{fig:case2Q}, respectively. 

\subsection{The oscillating droplet}

\begin{figure}[!htp]
\centering
\includegraphics[width=0.95\textwidth]{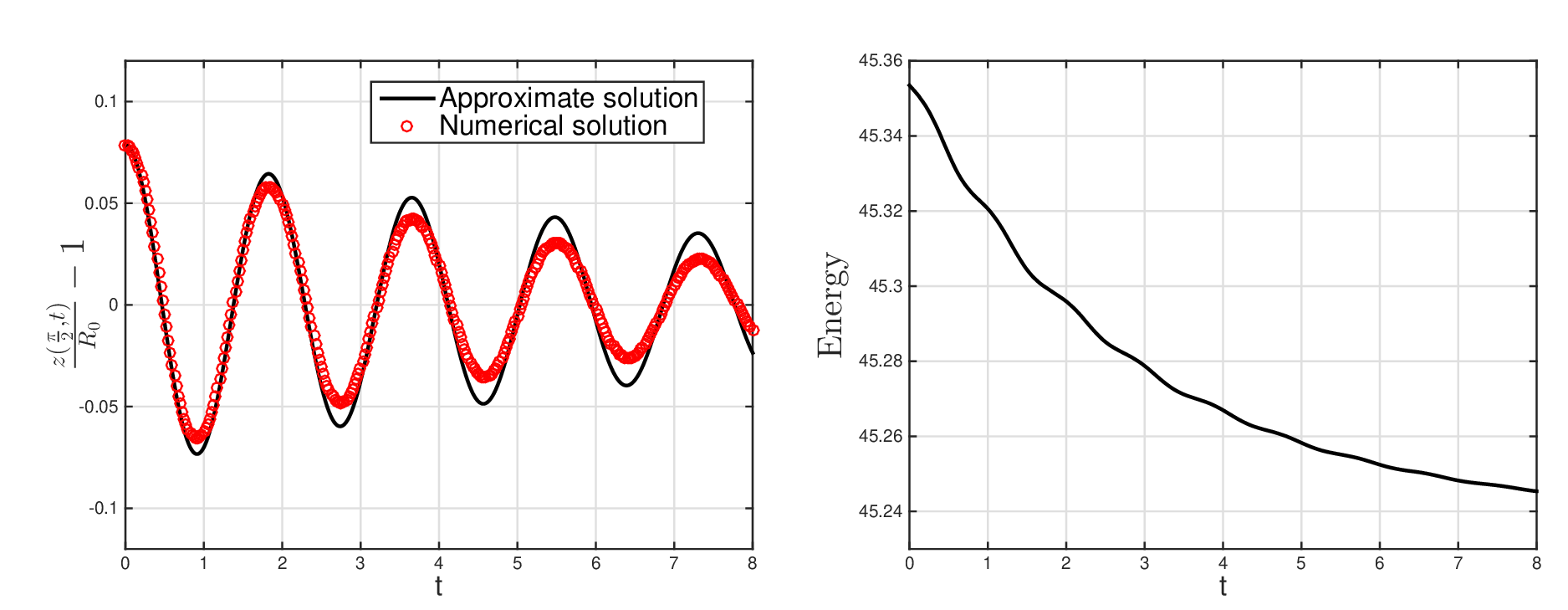}
\caption{The time evolution of the displacement of the upper endpoint of the generating curve on the $z$-axis (left panel) and the total energy (right panel), where we use $\EquiL^h$ with $2{\rm adapt}_{9,4}$.}
\label{fig:osd1}
\end{figure}

\begin{figure}[!htp]
\centering
\includegraphics[width=0.31\textwidth]{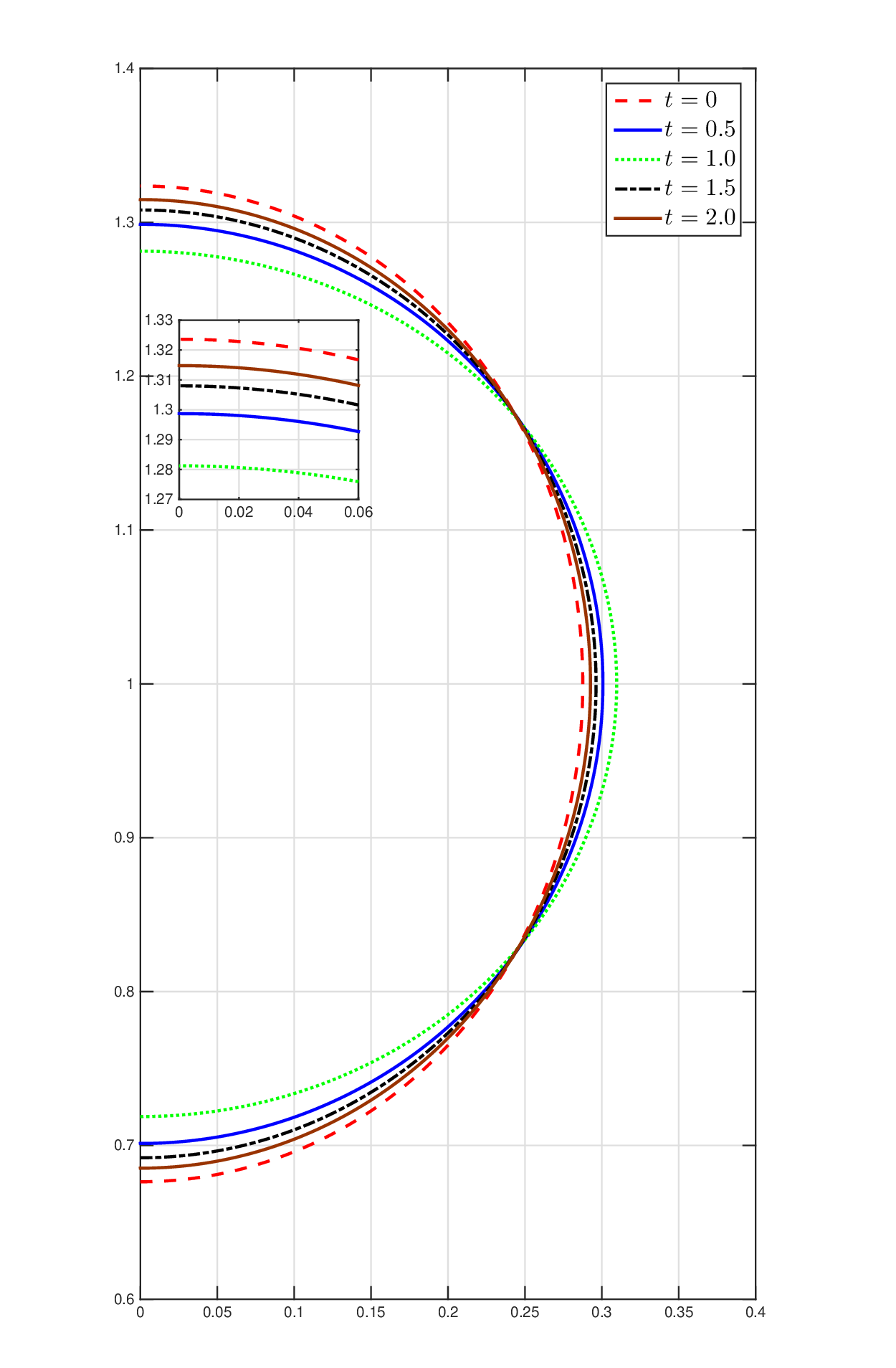}
\includegraphics[width=0.31\textwidth]{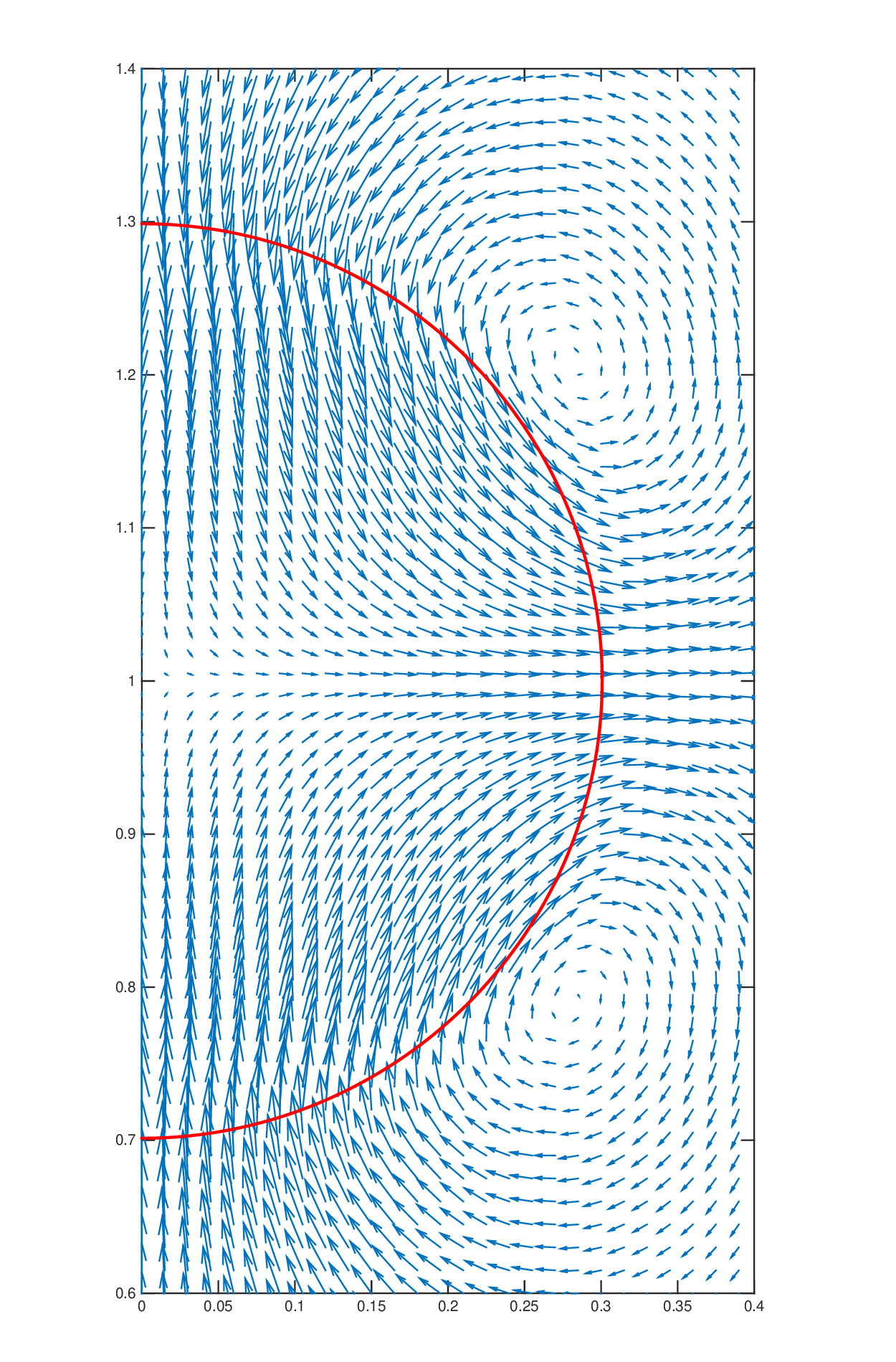}
\includegraphics[width=0.31\textwidth]{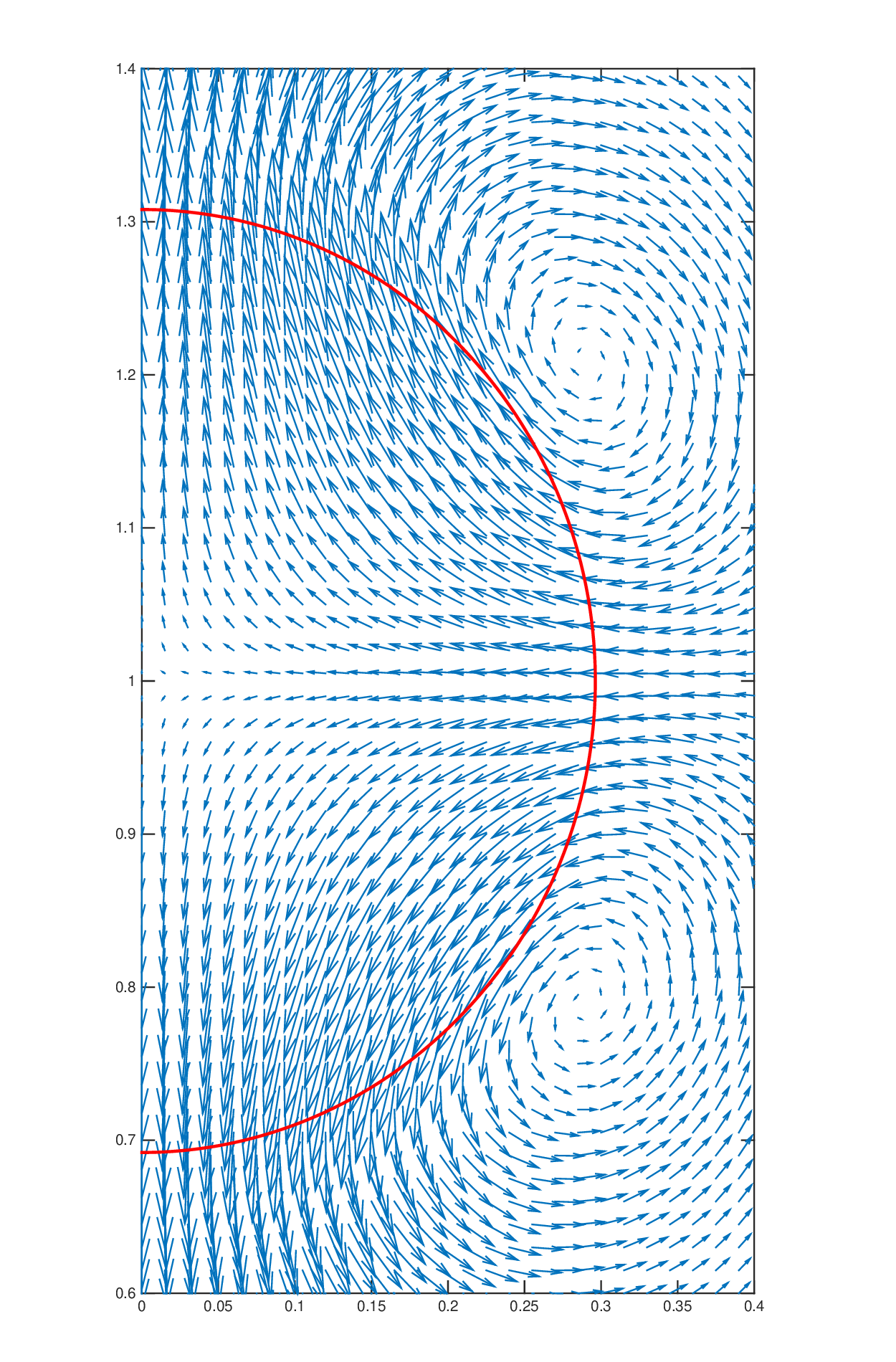}
\caption{Snapshots of the generating curves of the fluid interface at several times (left panel), the velocity fields at time $t=0.5$ (middle panel) and at time $t=1.5$ (right panel) for the experiment in Fig.~\ref{fig:osd1}.}
\label{fig:osd2}
\end{figure}

In this experiment, we study numerically the oscillation of a levitated droplet which is surrounded by a low density fluid. This is a classical problem and has been studied in detail in the literature, see \cite{Rayleigh1879capillary, Lamb81, Velentine1965,aalilija20}. Here we consider the oscillation of a 3d axisymmetric droplet, with the generating curve of the initial interface of the droplet given by
 \begin{equation}
\left\{\begin{array}{ll}r(\theta, 0) &= R_0\left[1 + \frac{\varepsilon_0}{2}(3\sin^2\theta-1)- \frac{1}{5}\varepsilon_0^2\right]\cos(\theta),\\
z(\theta, 0) &= R_0\left[1 + \frac{\varepsilon_0}{2}(3\sin^2\theta-1)- \frac{1}{5}\varepsilon_0^2\right]\sin(\theta) + 1.0,\nn
\end{array}\right.\quad \theta\in[-\tfrac{\pi}{2},\tfrac{\pi}{2}],
\end{equation}
where $R_0=0.3$ and $\varepsilon_0 = 0.08$.  For simplicity, we choose the physical parameters as
\begin{equation}
\rho_+=1, \quad\rho_- = 1000, \quad\mu_+=0.01,\quad\mu_-=2,\quad \gamma = 40,\quad \vec g = \vec 0.\nonumber
\end{equation}
The computational domain is $\mR = [0, 0.5]\times[0,2]$ with $\partial_1\mR =[0,0.5]\times\{0,2\}$ and $\partial_2\mR=\{0.5\}\times[0,2]$.
Then, on recalling \cite[(15b) and (38)]{aalilija20}, we note that the dynamic interface of the droplet can be approximated by
\begin{equation}
\left\{\begin{array}{ll}r(\theta, t) &= R_0\left[1 + \frac{\varepsilon(t)}{2}(3\sin^2\theta-1)- \frac{1}{5}\varepsilon(t)^2\right]\cos(\theta),\\
z(\theta, t) &= R_0\left[1 + \frac{\varepsilon(t)}{2}(3\sin^2\theta-1)- \frac{1}{5}\varepsilon(t)^2\right]\sin(\theta) + 1.0,
\end{array}\right.\quad \theta\in[-\tfrac{\pi}{2},\tfrac{\pi}{2}],\label{eq:otf}
\end{equation}
and $\varepsilon(t)$ is given by
\begin{equation}
\varepsilon(t)\approx \varepsilon_0\exp(-\lambda_{2} t)\cos(\omega_{2} t)\quad\mbox{with}\quad  \lambda_2 = \frac{5\mu_-}{\rho_-R_0^2}=\tfrac{1}{9},\quad \omega_2=\sqrt{\frac{8\gamma}{\rho_-R_0^3}}=\tfrac89 \sqrt{15}.\nn
\end{equation}
In Fig.~\ref{fig:osd1} our obtained numerical results are compared with the approximate solution in \eqref{eq:otf} for the displacement of the upper endpoints of the generating curve. Here a very good agreement is observed, which demonstrates the accuracy of our introduced methods. Finally, we  visualise the interface and the velocity fields in Fig.~\ref{fig:osd2}. 

\section{Conclusions}\label{sec:con}
We proposed four unfitted finite element methods for two-phase incompressible flow in an axisymmetric setting. The proposed methods combine an unfitted finite element method for the Navier--Stokes equation in the 2d-meridian halfplane together with a parametric finite element method in terms of the generating curve of the axisymmetric interface. We considered two different weak formulations with two possible approximations of the surface tension force. With suitable discretizations, we obtained an unconditionally stable method and a linear discretized method with the property of equidistribution. Furthermore, on utilizing suitable time-weighted discrete normals, we adapted the two introduced methods to achieve an exact volume preservation of the two phases.  Numerical results were presented to show the robustness of the introduced methods and to verify these good properties.

\section*{Funding} The work of Zhao was supported by the Alexander von Humboldt Foundation.

\section*{Data availibility}The datasets generated and/or analysed during the current study are available from the corresponding author on reasonable request.
\section*{Declarations}The authors certify that they have no affiliations with or involvement in any organization or entity with any financial interest or non-financial interest in the subject matter or materials discussed in this manuscript.

\appendix
\numberwithin{equation}{section}
\setcounter{equation}{0}
\section{Derivation of \eqref{eqn:weakform}}
\label{app:asywf}

Let $\vec x= (x,y,z)$ be the Cartesian coordinates, and $(r,z,\varphi)$ be the cylindrical coordinates such that
\begin{equation}
x = r\cos\varphi,\qquad y = r\sin\varphi,\nn
\end{equation}
where $r=\sqrt{x^2+y^2}$ is the radial distance and $\varphi$ is the  azimuthal angle. This gives rise to the  basis vectors $\big\{\vec e_r,\vec e_\varphi,\vec e_z\big\}$ 
\begin{equation}
\vec e_r = (\cos\varphi,\sin\varphi, 0)^T,\qquad \vec e_\varphi = (-\sin\varphi,\cos\varphi, 0)^T,\qquad\vec e_z=(0,0,1)^T,\label{eq:cbasis}
\end{equation}
which satisfy
\begin{equation}
\pp{\varphi}{\vec e_r} = \vec e_\varphi,\qquad \pp{\varphi}{\vec e_\varphi} = -\vec e_r,\qquad \pp{\varphi}{\vec e_z} = \vec 0,\qquad\pp{r}{\vec e_r} = \pp{r}{\vec e_\varphi}=\pp{r}{\vec e_z}=\vec 0.
\label{eq:rela}
\end{equation}
For any scalar function $\breve{b}(\vec x)$, we write in cylindrical coordinates $b(r,\varphi,z) = \breve{b}(\vec x)$. Then we have
\begin{equation}
\nabla \breve{b}(\vec x) = \pp{r}{b}\,\vec e_r + \frac{1}{r}\pp{\varphi}{b}\,\vec e_\varphi + \pp{z}{b}\,\vec e_z.
\end{equation}
For any vector-valued function $\breve{\vec b}(\vec x, t)$, similarly we write $\breve{\vec b}(\vec x)=\vec b(r,\varphi, z) = b^r\,\vec e_r + b^\varphi\,\vec e_\varphi + b^z\,\vec e_z$. Then using the relations in \eqref{eq:rela} gives
\begin{equation}
\nabla \breve{\vec b} =[\vec e_r\quad \vec e_\varphi \quad \vec e_z]\left[\begin{array}{lll}
\pp{r}{b^r} &\frac{1}{r}\pp{\varphi}{b^r}-\frac{b^\varphi}{r} &\pp{z}{b^r}\\[0.4em]
\pp{r}{b^\varphi}&\frac{b^r}{r}+ \frac{1}{r}\pp{\varphi}{b^\varphi} &\pp{z}{b^\varphi}\\[0.4em]
\pp{r}{b^z}&\frac{1}{r}\pp{\varphi}{b^z} &\pp{z}{b^z} 
\end{array}\right]
\left[\begin{array}{l}
(\vec e_r)^T\\[0.4em]
(\vec e_\varphi)^T\\[0.4em]
(\vec e_z)^T
\end{array}
\right].\label{eq:gradf}
\end{equation}
Moreover, for the divergence of $\breve{\vec b}$ it holds that
\begin{align}
\nabla\cdot\breve{\vec b} = \pp{r}{\vec b}\cdot\vec e_r + \frac{1}{r}\pp{\varphi}{\vec b}\cdot\vec e_\varphi + \pp{z}{\vec b}\cdot\vec e_z = \pp{r}{b^r} + \frac{b^r}{r}+\frac{1}{r}\pp{\varphi}{b^\varphi} +\pp{z}{b^z}. 
\label{eq:divf}
\end{align}

 We now follow the axisymmetric setting in \S\ref{sec:asyset} and denote 
 \[\breve{\vec u}(\vec x, t)=u^r(r,z,t)\,\vec e_r + u^z(r,z,t)\,\vec e_z,\qquad\breve{\vec\chi}(\vec x, t) =\chi^r(r,z,t)\,\vec e_r + \chi^z(r,z,t)\,\vec e_z.\] Then it follows from \eqref{eq:rhoc} that
\begin{align}
\Bigl(\rho\breve{\vec u},~\breve{\vec\chi}\Bigr)_\Omega &= 2\pi\int_{\mR_-(t)}\rho_-(u^r\,\chi^r + u^z\,\chi^z)r\,\drz +  2\pi\int_{\mR_+(t)}\rho_+(u^r\,\chi^r + u^z\,\chi^z)r\,\drz\nn\\
&=2\pi\int_{\mR}\rho_{_c}(u^r\,\chi^r + u^z\,\chi^z)r\,\drz=2\pi\,\Bigl(\rho_{_c}\vec u,~\vec\chi\,r\Bigr),
\end{align}
where $\vec u = (u^r, u^z)^T$ and $\vec\chi = (\chi^r,\chi^z)^T$ are vectors in the 2d-meridian halfplane $\mR$. Similarly we have
\begin{equation}
\Bigl(\rho\partial_t\breve{\vec u},~\breve{\vec\chi}\Bigr)_\Omega = 2\pi\,\Bigl(\rho_{_c}\partial_t\vec u,~\vec\chi\,r\bigr),\qquad \Bigl(\rho\breve{\vec u},~\partial_t\breve{\vec\chi}\Bigr)_\Omega = 2\pi\,\Bigl(\rho_{_c}\vec u,~\partial_t\vec\chi\,r\Bigr).
\end{equation}

On recalling \eqref{eq:divf} and \eqref{eq:nbc}, it holds that
\begin{equation}
\Bigl(\breve{q},\,\nabla\cdot\breve{\vec u}\Bigr) = 2\pi\int_\mR\left(\frac{\partial u^r}{\partial r} + \frac{u^r}{r} + \frac{\partial u^z}{\partial z}\right)q\,r\drz = 2\pi\,\Bigl(\bG\cdot([r\vec u],~q\Bigr).
\end{equation}
Using \eqref{eq:gradf} for $\breve{\vec u}$ and $\breve{\vec\chi}$, it is not difficult to show that
\begin{align}
\int_{\Omega_\pm(t)}\left([\breve{\vec u}\cdot\nabla]\breve{\vec u}\cdot\breve{\vec\chi} -[\breve{\vec u}\cdot\nabla]\breve{\vec\chi}\cdot\breve{\vec u}\right)\,\dV= 2\pi\int_{\mR_\pm(t)}\,\bigl([\vec u\cdot\bG]\vec u\cdot\vec\chi - [\vec u\cdot\bG]\vec\chi\cdot\vec u\bigr)\,r\,\drz.
\label{eq:axyi2}
\end{align}
Multiplying \eqref{eq:axyi2} with $\rho_\pm$ and combing these two equations gives rise to 
\begin{equation}
\Bigl(\rho,~[\breve{\vec u}\cdot\nabla]\breve{\vec u}\cdot\breve{\vec\chi} -[\breve{\vec u}\cdot\nabla]\breve{\vec\chi}\cdot\breve{\vec u}\Bigr)_\Omega =2\pi\,  \Bigl(\rho_{_c}r,~[\vec u\cdot\bG]\vec u\cdot\vec\chi - [\vec u\cdot\bG]\vec\chi\cdot\vec u\Bigr).
\end{equation}
Similarly we can compute
\begin{align}
\mat\tD(\breve{\vec u}) = \frac{1}{2}\left(\nabla\breve{\vec u}+ (\nabla\breve{\vec u})^T\right)=[\vec e_r\quad \vec e_\varphi \quad \vec e_z]\left[\begin{array}{lll}
\pp{r}{u^r} &0 &\frac{1}{2}\left(\pp{z}{u^r}+\pp{r}{u^z}\right)\\[0.4em]
0&\frac{u^r}{r} &0\\[0.4em]
\frac{1}{2}\left(\pp{z}{u^r}+\pp{r}{u^z}\right) &0 &\pp{z}{u^z} 
\end{array}\right]
\left[\begin{array}{l}
(\vec e_r)^T\\[0.4em]
(\vec e_\varphi)^T\\[0.4em]
(\vec e_z)^T
\end{array}
\right].\nn
\end{align}
This implies that
\begin{align}
&
\int_{\Omega_\pm(t)}\mat{\tD}(\breve{\vec u}):\mat{\tD}(\breve{\vec\chi})\,\dV\nn\\&\quad=2\pi\int_{\mR_\pm(t)}\left\{\pp{r}{u^r}\pp{r}{\chi^r} + \frac{1}{2}\left[\pp{z}{u^r} + \pp{r}{u^z}\right]\left[\pp{z}{\chi^r} + \pp{r}{\chi^z}\right]+ \,\pp{z}{u^z}\pp{z}{\chi^z} + \frac{u^r\,\chi^r}{r^2}\right\} r\,\drz\nn\\
&\quad=2\pi\int_{\mR_\pm(t)}\bD(\vec u):\bD(\vec\chi)\,r\drz +2\pi\int_{\mR_\pm(t)}(\vec u\cdot\vec e_1)(\vec\chi\cdot\vec e_1)\,r^{-1}\,\drz.\label{eq:axyi3}
\end{align}
Multiplying \eqref{eq:axyi3} with $\mu_\pm$ and combing these two equations yields that
\begin{equation}
\Bigl(\mu\,\mat{\tD}(\breve{\vec u}),~\mat{\tD}(\breve{\vec\chi})\Bigr)_\Omega = 2\pi\,\Bigl(\mu_{_c}r\,\bD(\vec u),~\bD(\vec\chi)\Bigr) +2\pi\, \Bigl(\mu_{_c}r^{-1}(\vec u\cdot\vec e_1),~(\vec\chi\cdot\vec e_1)\Bigr).
\end{equation}

On the axisymmetric surface $\mS(t)$, by \eqref{eq:mSp} we can set
\begin{equation}
\mathcal{\vv V}(\vec x, t) = \mZ^r_t\,\vec e_r + \mZ^z_t\,\vec e_z,\qquad\vec n_{_\mS}(\vec x, t) = -\mZ^z_s\,\vec e_r + \mZ^r_s\,\vec e_z,\qquad \breve{\vec\eta}(\vec x) = \eta^r\vec e_r + \eta^z\vec e_z\quad\mbox{for}\quad\vec x\in\mS(t).
\end{equation}
Denote $\zeta(\alpha, t) = \breve{\zeta}(\vec x, t)$ for $\vec x\in\mS(t)$, $\alpha\in\bI$. Then it is easy to get
\begin{equation}
\int_{\mS(t)}(\mathcal{\vv V}-\breve{\vec u})\cdot\vec n_{_\mS}\,\breve{\zeta}\,\rd S = 2\pi\int_{\Gamma(t)}(\vec\mZ_t - \vec u)\cdot\vec\nu\,\zeta\,(\vec\mZ\cdot\vec e_1)\rd s = 2\pi\,\Big\langle(\vec\mZ\cdot\vec e_1)(\vec\mZ_t-\vec u)\cdot\vec\nu,\,\zeta\,|\vec\mZ_\alpha|\Big\rangle.
\end{equation}
Similarly, it is straightforward to obtain
\begin{equation}
\int_{\mS(t)}\mathcal{H}\,\vec n_{_\mS}\cdot\breve{\vec\eta}\,\rd S =2\pi\int_{\Gamma(t)}\varkappa\,\vec\nu\cdot\vec\eta\;(\vec\mZ\cdot\vec e_1)\rd s = 2\pi\,\Big\langle(\vec\mZ\cdot\vec e_1)\,\varkappa\,\vec\nu,~\vec\eta\,|\vec\mZ_\alpha|\Big\rangle,
\end{equation}
where $\vec\eta=(\eta^r,~\eta^z)^T$ is the vector in the 2d-meridian plane. Finally, we have
\begin{align}
\int_{\mS(t)}\nabla_s\vec\id:\nabla_s\breve{\vec\eta}\,\rd S &= \int_{\mS(t)}\nabla_s\cdot\breve{\vec\eta}\,\rd S = 2\pi\int_{\Gamma(t)}(\vec\mZ\cdot\vec e_1)\left(\vec\mZ_s\cdot\vec\eta_s + \frac{\vec\eta\cdot\vec e_1}{\vec\mZ\cdot\vec e_1}\right)\rd s\nn\\
& = 2\pi \int_{\bI}(\vec\mZ\cdot\vec e_1)\vec\mZ_s\cdot\vec\eta_s\,|\vec\mZ_\alpha| + (\vec\eta\cdot\vec e_1)\,|\vec\mZ_\alpha|\,\rd\alpha,
\end{align}
where we used the fact 
\begin{equation}
\nabla_s\cdot\breve{\vec\eta} = \mZ^r_s\eta^r_s + \mZ^z_s\eta^z_s + \frac{\eta^r}{\vec\mZ\cdot\vec e_1}.
\end{equation}
Collecting the above results, we obtain \eqref{eqn:weakform} from \eqref{eqn:3dwf}.

\bibliographystyle{model1b-num-names}
\bibliography{bib}
\end{document}